\documentclass[a4paper,reqno,10pt]{amsart}

\usepackage{footmisc}

\usepackage{amsmath}
\usepackage{amsthm}
\usepackage{amssymb}
\usepackage{amscd} 
\allowdisplaybreaks

\usepackage{bbm} 
\usepackage[scaled=.98]{BOONDOX-uprscr}
\usepackage[OMLmathsfit]{isomath}
\usepackage{mathabx}
\usepackage{txfonts}

\usepackage[T1]{fontenc} 
\usepackage{lmodern}

\usepackage{graphicx} 

\usepackage{imakeidx}
\makeindex[name=symbol,title={Symbol Index}] 
\makeindex[name=author,title={Author Index}] 
\makeindex[name=subject,title={Subject Index}] 

\usepackage{verbatim} 

\usepackage[Lenny]{fncychap}

\usepackage{sqrcaps} 

\usepackage{url}
\usepackage{comment}
\usepackage[svgnames]{xcolor}
\usepackage{pdfcolmk}

\usepackage{tikz}
\usepackage{pifont}

\usepackage{enumitem}
\setlist{noitemsep}

\usepackage{stackengine}

\usepackage{pgfplots}

\usepackage{hyphenat}

\swapnumbers 

\newtheoremstyle{exercise} 
  {3pt} 
  {3pt} 
  {\small\rmfamily} 
  {
} 
  {\rmfamily\scshape} 
  {.} 
  {.5em} 
  {} 

\newtheoremstyle{newplain}
  {5pt}
  {5pt}
  {\itshape}
  {}
  {\rmfamily\scshape}
  { ---}
  {.5em}
  {}

\newtheoremstyle{newremark}
  {5pt}
  {5pt}
  {\rmfamily}
  {}
  {\rmfamily\scshape}
  { ---}
  {.5em}
  {}

\theoremstyle{newplain}
\newtheorem*{Theorem*}{Theorem} 

\theoremstyle{newplain}
\newtheorem{Theorem}{Theorem}

\newtheorem{Corollary}[Theorem]{Corollary}
\newtheorem{Proposition}[Theorem]{Proposition}
\newtheorem{Conjecture}[Theorem]{Conjecture}
\newtheorem{Definition}[Theorem]{Definition}

\theoremstyle{newremark}
\newtheorem{Empty}[Theorem]{}
\newtheorem{Remark}[Theorem]{Remark}

\newtheorem{Claim}[Theorem]{Claim}

\theoremstyle{exercise}

\numberwithin{Theorem}{section}
\numberwithin{equation}{section}

\newcommand{\N}{\mathbb{N}} 

\newcommand{\R}{\mathbb{R}} 
\newcommand{\Rm}{\R^m}
\newcommand{\Rn}{\R^n}
\newcommand{\Rp}{\R^p}
\newcommand{\Z}{\mathbb{Z}} 

\newcommand{\ind}{\mathbbm{1}} 

\newcommand{\calD}{\mathscr{D}}
\newcommand{\calE}{\mathscr{E}}
\newcommand{\calF}{\mathscr{F}}
\newcommand{\calG}{\mathscr{G}}

\newcommand{\calL}{\mathscr{L}}

\newcommand{\calP}{\mathscr{P}}
\newcommand{\calQ}{\mathscr{Q}}

\newcommand{\calT}{\mathscr{T}}

\newcommand{\bkappa}{\boldsymbol{\kappa}}

\newcommand{\bmu}{\boldsymbol{\mu}}

\newcommand{\btau}{\boldsymbol{\tau}}

\newcommand{\bB}{\mathbf{B}}
\newcommand{\bC}{\mathbf{C}}

\newcommand{\bF}{\mathbf{F}}
\newcommand{\bG}{\mathbf{G}}

\newcommand{\bL}{\mathbf{L}}
\newcommand{\bM}{\mathbf{M}}
\newcommand{\bN}{\mathbf{N}}

\newcommand{\bP}{\mathbf{P}}

\newcommand{\bc}{\pmb{c}}

\DeclareMathOperator{\rmcard}{\mathrm{card}} 
\DeclareMathOperator{\rmclos}{\mathrm{clos}} 
\DeclareMathOperator{\rmcpt}{\mathrm{cpt}} 
\DeclareMathOperator{\rmdiam}{\mathrm{diam}} 

\DeclareMathOperator{\rmdist}{\mathrm{dist}} 
\DeclareMathOperator{\rmdiv}{\mathrm{div}} 

\DeclareMathOperator{\rmid}{\mathrm{id}} 
\DeclareMathOperator{\rmim}{\mathrm{im}} 
\DeclareMathOperator{\rmint}{\mathrm{int}} 

\DeclareMathOperator{\rmLip}{\mathrm{Lip}} 

\DeclareMathOperator{\rmosc}{\mathrm{osc}} 
\DeclareMathOperator{\rmspt}{\mathrm{spt}} 

\newcommand{\ip}{\,\,\begin{picture}(-1,1)(-1,-2.5)\circle*{2}\end{picture}\;\,\,} 

\newcommand{\hel} {
\hskip2.5pt{\vrule height7pt width.5pt depth0pt}
\hskip-.2pt\vbox{\hrule height.5pt width7pt depth0pt}
\, }

\newcommand{\cone} {
\times
\mspace{-15.0mu} 
\times 
\medspace }

\newcommand{\lseg}{\boldsymbol{[}\!\boldsymbol{[}}
\newcommand{\rseg}{\boldsymbol{]}\!\boldsymbol{]}}

\newcommand{\cqfd} {
\renewcommand{\qedsymbol}{$\blacksquare$}
\qed
\renewcommand{\qedsymbol}{$\square$} }

\newcommand{\bin}[2]{
\begin{pmatrix} #1 \\
#2 
\end{pmatrix}}

\def\XXint#1#2#3{{%
\setbox0=\hbox{$#1{#2#3}{\int}$}
\vcenter{\hbox{$#2#3$}}\kern-.5\wd0}}

\renewcommand{\em}{\bf}

\newcommand{\veps}{\varepsilon}

\newcommand{\vphi}{\varphi}

\newcommand{\la}{\langle}
\newcommand{\ra}{\rangle}

\newcommand{\ie}{{\it i.e.}\ }
\newcommand{\eg}{{\it e.g.}\ }

\renewcommand{\leq}{\leqslant}
\renewcommand{\geq}{\geqslant}
\renewcommand{\subset}{\subseteq}

\newlength{\drop}

\usepackage{trajan}

\begin{document}



\title[Lipschitz-free spaces]{Quantified compactness\\ in Lipschitz-free spaces of $[-1,1]^n$}

\author[Th. De Pauw]{Thierry De Pauw}

\address{Institute for Theoretical Sciences / School of Science, Westlake University\\
No. 600, Dunyu Road, Xihu District, Hangzhou, Zhejiang, 310030, China}

\email{thierry.depauw@westlake.edu.cn}





\begin{abstract}
We show that the members of the Lipschitz-free space of $[-1,1]^n$ are exactly the 0-dimensional flat currents whose ``boundary'' vanishes.
The connection with normal and flat currents allows to use the Federer-Fleming compactness and deformation theorems in this context.
We characterize the compact subsets of this Lipschitz-free space and we quantify their $\veps$-entropy.
\end{abstract}

\maketitle

\tableofcontents


\section{Introduction}

Given a metric space $X$ and $o \in X$ we let $\rmLip_o(X)$ be the Banach space consisting of those Lipschitz functions $u : X \to \R$ such that $u(o)=0$ equipped with its norm $\|u\|_L = \rmLip u$.
The map $X \to \rmLip_o(X)^* : x \mapsto \delta_x$ is an isometric embedding.
The closure of the vector subspace of $\rmLip_o(X)^*$ spanned by $\{ \delta_x : x \in X \}$ is called the Lipschitz-free space of $X$ and denoted $\calF(X)$, \ref{P.LF}.
Lipschitz-free spaces have been studied for several decades, see \eg \cite{GOD.KAL.03,GOD.15,WEAVER.2nd}, from the point of view of Banach space theory.
Here, we concentrate on studying {\it members} of $\calF(X)$ and we introduce tools from geometric measure theory to this end.
\par 
The construction of $\calF(X)$ is based on linear combinations of $\delta_x$, $x \in X$, but it is misleading to think of those as measures.
Indeed, the completion process defining $\calF(X)$ is not understood with respect to mass norm (duality with $C_0(X)$ in case $X$ is locally compact) but rather with respect to the norm $\|\cdot\|_L^*$ (duality with $\rmLip_o(X)$).
This has been noted for a long time as well as the connection with, \eg, the Wasserstein distance and dually with mass transportation.
Since we will refer to the deformation theorem, we focus on the case $X = [-1,1]^n \subset \Rn$ (and $o=0$) rather than more general metric spaces.
Note that $\calF(\Rn)$ and $\calF([-1,1]^n)$ are isomorphic Banach spaces, according to \cite{KAU.15}.
\par 
Let $E \subset \calF([-1,1]^n)$ be the vector subspace spanned by the image of the isometric embedding $[-1,1]^n \to \rmLip_o([-1,1]^n)$.
With each $\alpha = \sum_i \theta_i \cdot \delta_{x_i} \in E$ we associate a distribution $\iota(\alpha) = \sum_i \theta_i \cdot \lseg x_i \rseg - \left( \sum_i \theta_i \right) \cdot \lseg 0 \rseg$, where $\la \vphi , \lseg x \rseg \ra  = \vphi(x)$ for all test functions $\vphi$.
The distribution $T = \iota(\alpha)$ satisfies the following properties:
\begin{equation*}
\|\alpha\|_L^* = \bG(T),
\end{equation*}
where 
\begin{equation}
\label{eq.intro.1}
\bG(T) = \sup \{ \la \vphi , T \ra : \vphi \text{ is a test function and } \bM(d\vphi) \leq 1 \},
\end{equation}
($\bM(d\vphi)$ is simply the maximum norm of $d\vphi$, hence, equals $\rmLip \vphi$, see \ref{P.M}) and 
\begin{equation*}
\chi(T) = 0,
\end{equation*}
where $\chi(T) = \la \vphi , T \ra$ for any test function $\vphi$ such that $\rmspt T \subset \rmint \{ \vphi = 1 \}$ (such $\vphi$ exists, since $T$ has compact support).
This naturally leads to consider the space of distributions arising as $\iota(\alpha)$, for some $\alpha \in E$, with respect to the norm $\bG$ and to identify the limiting distributions.
\par 
Flat currents (also called flat chains) have been introduced by H. Whitney in \cite{WHITNEY}, see \ref{P.FNFC}.
Viewing members of $\calF([-1,1]^n)$ as mere distributions allows one to push them forward only by smooth mappings, though, from their very nature they can be pushed forward by Lipschitz mappings.
Flat currents, too.
In fact, flat currents of varying dimensions come equipped with a structure of chain complex (corresponding to their boundary operator) and morphisms thereof (corresponding to push-forward by Lipschitz mappings).
Furthermore, flat currents do not need to be of finite mass (\ie to be multi-vector valued Borel measures).
A typical example of a 0-dimensional flat current\footnote{0-dimensional currents are distributions} which is also a $\bG$-limit of members of $\iota(E)$ is
\begin{equation*}
\partial \sum_{j=1}^\infty \lseg 0, j^{-2} \rseg.
\end{equation*}
Based on this, we abbreviate $\bC^n = [-1,1]^n$ and we define $\calG(\bC^n)$ to consist of those 0-dimensional flat currents $T$ supported in $\bC^n$ such that $\chi(T)=0$.
We call these the 0-dimensional flat cycles.
\par 
The fact that $\bG(T)$ defined in \eqref{eq.intro.1} equals the infimum of $\bM(S)$ over all 1-dimensional currents $S$ with $\partial S = T$ is an easy consequence of Hahn-Banach's theorem, \ref{3.1} and \ref{3.2}.
That, for those $T$ that arise as $\bG$-limits of members of $\iota(E)$, $\bG(T)$ be given by 
\begin{multline}
\label{eq.intro.2}
\bG(T) = \inf \{ \bM(S) : T = \partial S \text{ and $T$ is a 1-dimensional flat current}\\ \text{ supported in } \bC^n \},
\end{multline}
is not as straightforward, see \ref{5.2}(C).
This follows from representation theorems of 0-dimensional flat cycles in section \ref{sec.rep}.
One of these, specifically \ref{4.2}, is akin a representation theorem known for members of $\calF(\bC^n)$, see \cite{MAR.KAL.KAP.17,GOD.LER.18}.
These authors show that $\calF(\bC^n)$ is isometrically isomorphic to the space of distributions $\rmdiv v$, for some $v \in L_1(\bC^n;\Rn)$, equipped with the quotient norm.
The proof relies on the fact that both $\calF(\bC^n)$ and the aforementioned space of distributional divergences are a predual of $\rmLip_o(\bC^n)$ (in this respect, see also the appendix of \cite{MAL.14}) and the uniqueness of such predual.
Similarly, we show that each $T \in \calG(\bC^n)$ is of the type $T = \partial (\calL^n \wedge \eta)$ for some Lebesgue-summable $\eta : \bC^n \to \bigwedge_1 \Rn$, with near correspondence of the relevant norms, but our argument in the proof of \ref{4.2} is somewhat constructive. 
We are then able to conclude that $\calF(\bC^n)[\|\cdot\|_L^*]$ and $\calG(\bC^n)[\bG]$ are isometrically isomorphic Banach spaces.
\par 
Since $T \in \calG(\bC^n)$ may have infinite mass, we introduce a measure of how much it misses having finite mass.
Specifically, we let $\bkappa(T,\veps)$ be the infimum of $\bM(\hat{T})$ over all $\hat{T} \in \calG(\bC^n)$ so that $\bG(T - \hat{T}) < \veps$.
The elementary properties of $\veps \mapsto \bkappa(T,\veps)$ are studied in \ref{8.1} and \ref{8.0}.
Referring to the Federer-Fleming compactness theorem of normal (1-dimensional) currents, in \ref{8.3} we characterize the relatively $\bG$-compact subsets of $\calG(\bC^n)$ (among those that are $\bG$-bounded) by means of $\bkappa(T,\veps)$, akin the Arzel\`a-Ascoli theorem for $C(X)$ or the Riesz-Kolmogoroff theorem for $\bL_1(\Rn)$.
Thus, $\bkappa(T,\veps)$ plays a similar role to the modulus of continuity (in Arzel\`a-Ascoli) or the modulus of smoothness in $\bL_1$ (in Riesz-Kolmogoroff). 
We also show in \ref{8.4} how $\bkappa(T,\veps)$ is related to the modulus of smoothness in $\bL_1$ of $\eta$, where $T = \partial(\calL^n \wedge \eta)$.
\par 
Most of this could be extended to the realm of a general complete metric space (rather than merely the unit cube of $\Rn$) by means of the general tools of geometric measure theory set forth in the past decades, namely metric currents \cite{AMB.KIR.00}, local metric currents \cite{LAN.11}, flat chains in Banach spaces \cite{ADA.08}, or flat chains in metric spaces \cite{DEP.HAR.07}.
However, in the last two sections of this paper we {\it quantify} the $\veps$-entropy of compact subsets of $\calG(\bC^n)$ in terms of an upper bound for $\bG$, a control on $\bkappa(\cdot,\veps)$, and the dimension $n$, see definition \ref{8.5} and \ref{10.4}.
This quantitative information is obtained by means of the {\it deformation theorem} valid only in $\Rn$ (though, see \cite[A.1 and A.2]{DEP.14}), which motivates our desire to restrict, here, to the classical theory of Federer-Fleming.
Regarding the deformation theorem, the reader not yet acquainted with it is invited to consult \cite{SIMON}, \cite{MORGAN}, \cite{WHI.99.deformation}, or the original \cite[4.2.9]{GMT}.
\par 
Our main quantitative result \ref{10.4} is as follows.
Assume that:
\begin{enumerate}
\item $0 < \veps \leq 1$;
\item $\calT \subset \calG(\bC^n)$;
\item $\Gamma > 0$ and $\bG(T) \leq \Gamma$ for all $T \in \calT$;
\item $k$ is a positive integer such that
\begin{equation*}
 \bc(n) \cdot \left( \Gamma + \bkappa \left( T, \frac{\veps}{6} \right) \right) < k \cdot \left( \frac{\veps}{6} \right)
\end{equation*}
for all $T \in \calT$;
\item $\calT$ is an $\veps$-net, \ie
\begin{equation*}
(\forall T_1 \in \calT)(\forall T_2 \in \calT): T_1 \neq T_2 \Rightarrow \bG(T_1-T_2) \geq \veps.
\end{equation*}
\end{enumerate}
Then
\begin{equation*}
\rmcard \calT \leq [22 \cdot n \cdot k ]^{(2k+1)^n}.
\end{equation*}
\par 
The notations $\bG$ and $\calG(\bC^n)$ are chosen in reference to that used by F.J. Almgren (see \eg \cite{ALM.86}) in his study of multi-points and multiple-valued functions (probably because $\bG$ follows $\bF$ in the alphabet and the obvious link with flat norm).
Among many references about multiple-valued functions, see \eg \cite{BOU.DEP.GOB,DEL.SPA.11}.
\par 
In the next section we try to gather the necessary material borrowed to the theory of normal and flat currents.
Perhaps \ref{P.AM} is new (augmentation map) as well \ref{P.THM}.
We have made an effort to stick with the notations of H. Federer's treatise \cite{GMT} in an attempt to facilitate the reader's immersion in the latter if necessary.
\par 
It is a pleasure to record useful discussions with Ph. Bouafia and G. Godefroy during the preparation of this paper.

\section{Preliminaries}

\begin{Empty}
We let $\R^+$ be the set of non-negative real numbers and $\R^+_0$ be the set of positive real numbers.
Similarly, $\N_0$ is the set of positive integers.
By $E[\nu]$ we mean a pair $(E,\nu)$ where $E$ is a real vector space and $\nu$ is a norm on $E$.
In that case we let $E^*$ be the dual of $E[\nu]$ and $\nu^*$ its canonical norm.
We denote as $E \times E^* \to \R : (x,x^*) \mapsto \la x, x^* \ra$ the corresponding duality.
\end{Empty}

\begin{Empty}[$\Rn$, $\ip$, $\bC^n$, $\pi$, and $\calL^n$]
\label{P.ES}
Throughout, $n$ is a positive integer.
When $x,y \in \Rn$ we let $x \ip y$ be the canonical inner product of $x$ and $y$ and we let $|x| = \sqrt{x \ip x}$ be the Euclidean norm of $x$.
The main structure of a metric space we will consider on $\Rn$ is the one induced by the norm $|\cdot|$.
In particular, Lipschitz constants are {\it always} meant with respect to the Euclidean norm.
Euclidean closed balls are denoted $\bB(a,r)$.
On one occasion (the proof of \ref{4.2}) we will use the maximum norm $|\cdot|_\infty$ in $\Rn$ and the corresponding balls will be denoted $\bB_\infty(a,r)$.
Furthermore, we abbreviate $\bC^n = [-1,1]^n = \bB_\infty(0,1)$ and we let $\pi : \Rn \to \Rn$ be the nearest point projection on $\bC^n$.
Thus, $\rmLip \pi = 1$.
Finally, $\calL^n$ is the Lebesgue measure on $\Rn$.
\end{Empty}

\begin{Empty}[Lipschitz-free spaces]
\label{P.LF}
Let $(X,d)$ be a non-empty metric space and $o \in X$.
We let $\rmLip_o(X)$ be the vector space consisting of those Lipschitz functions $u : X \to \R$ such that $f(o)=0$.
Equipped with the norm $u \mapsto \|u\|_L = \rmLip u$, this is a Banach space.
Letting $\la u , \delta_x \ra = u(x)$, $x \in X$, defines a map $\delta : X \to \rmLip_o(X)[\|\cdot\|_L]^*$.
Clearly, $\|\delta_x - \delta_y\|_L^* \leq d(x,y)$.
Considering $u_y(x) = d(y,x)-d(y,o)$, we observe that $u_y \in \rmLip_o(X)$ and $\la u_y , \delta_x - \delta_y \ra = d(x,y)$.
Accordingly, $\delta$ is an isometry.
We define the {\em Lipschitz-free space of the pointed metric space $(X,o)$} to be the closure in $\rmLip_o(X)[\|\cdot\|_L]^*$ of the vector subspace generated by the image of $\delta$.
This is a Banach space and we denote it by $\calF(X)[\|\cdot\|_L^*]$ or simply $\calF(X)$.
In this paper, we study $\calF(\bC^n)$ (we choose $o \in \bC^n$ to be the origin of $\Rn$).
\end{Empty}

\begin{Empty}[Differential forms and currents]
\label{P.TFC}
In this preliminary section we shall recall some notions and results about differential forms of degree $m$ and $m$-dimensional currents.
Here, $m$ is a non-negative integer.
However, starting from next section we shall only consider the cases $m=0$ and $m=1$.
Thus, the reader who does not want to acquaint themself with general values of $m$ may safely restrict to considering only those two particular cases.
\par 
We let $\bigwedge_m \Rn$ be the real vector space consisting of $m$-vectors of $\Rn$ and $\bigwedge^m \Rn$ its dual.
Both are equipped with an inner product structure -- and corresponding norm -- induced by that of $\Rn$ in a canonical way.
In particular, $\bigwedge_0 \Rn = \R \cong \bigwedge^0 \Rn$ and $\bigwedge_1 \Rn = \Rn \cong \bigwedge^1 \Rn$.
A {\em differential form of degree $m$} is a $C^\infty$ map $\omega : \Rn \to \bigwedge^m \Rn$.
Its {\em support} $\rmspt \omega$ is the closure of $\Rn \cap \{ x : \omega(x) \neq 0 \}$.
The vector space $\calE^m(\Rn)$ consists of all differential forms of degree $m$ and its vector subspace $\calD^m(\Rn)$ consists of those that have compact support.
The members of $\calD^0(\Rn)$ are often called {\em test functions}.
We usually denote them by the letter $\vphi$ or $\psi$.
\par 
The space $\calD^m(\Rn)$ is equipped with a locally convex topology described \eg in \cite[Chap. 6]{RUDIN}.
Its topological dual is denoted $\calD_m(\Rn)$ and its members are called {\em $m$-dimensional currents}.
In particular, members of $\calD_0(\Rn)$, which are usually called {\em distributions}, are also called $0$-dimensional currents.
The corresponding duality is denoted $\calD^m(\Rn) \times \calD_m(\Rn) : (\omega,T) \mapsto \la \omega , T \ra$.
The support of $T \in \calD_m(\Rn)$, denoted $\rmspt T$, is the complement of the largest open set $U \subset \Rn$ such that $\la \omega , T \ra = 0$ whenever $\rmspt \omega \subset U$.
\end{Empty}

\begin{Empty}[Exterior derivative and pull-back]
\label{P.EDPB}
We let $d : \calD^m(\Rn) \to \calD^{m+1}(\Rn)$ be the continuous linear operator that sends $\omega$ to its {\em exterior derivative} $d\omega$.
In particular, if $m=0$ then $d\vphi(x) = \sum_{i=1}^n D_i\vphi(x) \cdot dx_i$, where $D_i\vphi(x)$ is the $i^{th}$ partial derivative of $\vphi$ at $x$ and $dx_i$ is the $i^{th}$ coordinate function in $\Rn$. 
Note that $\rmspt d\omega \subset \rmspt \omega$.
\par 
Let $f : \Rn \to \Rp$ be $C^\infty$ and let $f^\# : \calE^m(\Rp) \to \calE^m(\Rn)$ be the linear operator that sends $\omega$ to its {\em pull-back} $f^\# \omega$.
In particular, if $m=0$ then $f^\# \omega = \omega \circ f$.
Note that $\rmspt f^\# \omega \subset f^{-1}(\rmspt \omega)$.
We say that $f$ is {\em proper} if $f^{-1}(K)$ is compact whenever $K$ is compact.
If $f$ is smooth and proper then $f^\#$ maps $\calD^m(\Rp)$ to $\calD^m(\Rn)$ and is continuous.
In that case, $d \circ f^\# = f^\# \circ d$.
\end{Empty}

\begin{Empty}[Boundary and push-forward]
\label{P.BPF}
If $m \geq 1$ then we let $\partial : \calD_m(\Rn) \to \calD_{m-1}(\Rn)$ be adjoint to the corresponding exterior derivative operator.
We call $\partial T$ the {\em boundary} of $T$.
Note that $\rmspt \partial T \subset \rmspt T$.
An ersatz of boundary for distributions (the case $m=0$) will be defined in \ref{P.AM}.
\par 
If $f : \Rn \to \Rp$ is $C^\infty$ and proper then we let $f_\# : \calD_m(\Rn) \to \calD_m(\Rp)$ be adjoint to the corresponding pull-back operator.
We call $f_\#T$ the {\em push-forward of $T$ by $f$}.
Note that $\rmspt f_\# T \subset f(\rmspt T)$.
Furthermore, $\partial \circ f_\# = f_\# \circ \partial$.
\par 
If $f$ is $C^\infty$ but not proper we define $f_\#T$ as follows for those $T \in \calD^m(\Rn)$ that have {\it compact support}.
Letting $\psi$ be a test function such that $\rmspt T \subset \rmint \{ \psi = 1 \}$ we set $\la \omega , f_\# T \ra = \la \psi \cdot f^\# \omega , T \ra$ for all $\omega \in \calD^m(\Rp)$.
One first checks that the definition is independent of the choice of $\psi$ and then that $T \mapsto f_\# T$ is linear and commutes with $\partial$ and that $\rmspt f_\# T \subset f(\rmspt T)$.
\end{Empty}

\begin{Empty}[Mass]
\label{P.M}
The {\em comass}\footnote{In case $2 \leq m \leq n-2$ this differs from and is equivalent to the comass defined in \cite[1.8.1 and 4.1.7]{GMT}. The difference will play no r\^ole here, since we focus on $m=0$ and $m=1$.} of $\omega \in \calD^m(\Rn)$ is $\bM(\omega) = \max \{ |\omega(x)| : x \in \Rn \} < \infty$.
Notice that $\bM(d\vphi) = \rmLip \vphi$ for all $\vphi \in \calD^0(\Rn)$.
The {\em mass} of a linear functional $T : \calD^m(\Rn) \to \R$ is $\bM(T) = \sup \{ \la \omega , T \ra : \omega \in \calD^m(\Rn) \text{ and } \bM(\omega) \leq 1 \} \in [0,\infty]$.
If $\bM(T) < \infty$ then $T$ is an $m$-dimensional current but the converse fails.
If $T$ is a distribution then the condition that $T$ have finite mass is also referred to as $T$ having order 0.
\par 
It follows from the Riesz-Markoff representation theorem that $T \in \calD_m(\Rn)$ has finite mass if and only if there exist a finite Borel measure $\|T\|$ on $\Rn$ and a Borel-measurable unit multi-vector field $\vec{T} : \Rn \to \bigwedge_m\Rn$ such that $\la \omega , T \ra = \int_{\Rn} \la \vec{T} , \omega \ra d\|T\|$ for all $\omega \in \calD^m(\Rn)$.
Furthermore, $\|T\|$ is uniquely determined by $T$ and $\vec{T}$ is uniquely determined $\|T\|$-almost everywhere.
\par 
Observe that $\bM(f_\# T) \leq \|Df\|_\infty^m \bM(T)$ if either $f$ is smooth and proper or $f$ is smooth and $T$ has compact support.
\par 
In case $m=0$, $\vec{T}$ is $\{-1,1\}$-valued $\|T\|$-almost everywhere.
Therefore, a finite mass distribution $T$ corresponds to the Borel signed measure $\mu_T(B) = \int_{B} \vec{T} d\|T\|$ by means of the formula $\la \vphi , T \ra = \int_{\Rn} \vphi d\mu_T$.
\par 
Finally, it will be useful to note that the vector space $\calD^m(\Rn) \cap \{ T : \bM(T) < \infty \}$ equipped with the norm $\bM$ is complete and that $\bM : \calD_m(\Rn) \to [0,\infty]$ is sequentially weakly* lower semi-continuous.
\end{Empty}

\begin{Empty}[Normal currents]
\label{P.NC}
If $T \in \calD_m(\Rn)$ and $m \geq 1$ then $\partial T$ may not have finite mass even though $\bM(T) < \infty$.
For instance, upon letting $\lseg a,b \rseg \in \calD_1(\R)$ be defined by $\la \vphi \cdot dx , \lseg a,b \rseg \ra = \int_a^b \vphi d\calL^1$, we note that $T = \sum_{j=1}^\infty \lseg -j^{-2} , j^{-2} \rseg$ defines a 1-dimensional current of finite mass in $\R$. 
Let $\lseg a \rseg$ be the distribution in $\R$ defined by $\la \vphi , \lseg a \rseg \ra = \vphi(a)$ and likewise $\lseg b \rseg$.
Since $\partial \lseg a,b \rseg = \lseg b \rseg - \lseg a \rseg$, according to the fundamental theorem of calculus, the partial sums $\sum_{j=1}^k \lseg j^{-2} \rseg - \lseg -j^{-2} \rseg$ converge weakly* to $\partial T$ as $k \to \infty$.
It is easy to check that $\bM(\partial T) = \infty$.
\par 
For $T \in \calD_m(\Rn)$ we let $\bN(T) = \bM(T) + \bM(\partial T)$ if $m \geq 1$ and $\bN(T) = \bM(T)$ if $m=0$.
We say that $T$ is {\em normal}\footnote{This is consistent with H. Federer's definition. One could argue that the appropriate definition in this preliminary section would be to only require that $\bN(T) < \infty$ and then observe that, among these, those having compact support constitute an $\bF$-dense subspace, see \ref{P.FNFC} for the definition of $\bF$.} if $T$ has compact support and if $\bN(T) < \infty$ and we let $\bN_m(\Rn)$ be the space of normal $m$-dimensional currents in $\Rn$.
Thus, $\bN_0(\Rn)$ consists of those distributions in $\Rn$ associated with a compactly supported Borel signed measure.
The definition is designed so that $\partial T \in \bN_{m-1}(\Rn)$ whenever $T \in \bN_m(\Rn)$ and $m \geq 1$.
\end{Empty}

\begin{Empty}[The cone construction]
\label{P.CONE}
Let $T \in \bN_m(\Rn)$ and $a \in \Rn$. 
There exists a current $\lseg a \rseg \cone T \in \bN_{m+1}(\Rn)$ called {\em the cone with vertex $a$ and base $T$}.
It is defined by means of the affine homotopy from the constant map $a$ to the identity of $\Rn$, \cite[4.1.11 middle of p.365]{GMT}.
For our purposes, we note that 
\begin{equation*}
\rmspt (\lseg a \rseg \cone T) \subset \mathrm{convex} ( \{a\} \cup \rmspt T)
\end{equation*}
and
\begin{equation*}
\bM(\lseg a \rseg \cone T) \leq \rmdiam  ( \{a\} \cup \rmspt T) \cdot \bM(T)
\end{equation*}
(it is important that $T$ have compact support for $\lseg a \rseg \cone T$ to make sense), and, in case $m=0$,
\begin{equation*}
\partial (\lseg a \rseg \cone T) = T - \la \vphi , T \ra \cdot \lseg a \rseg,
\end{equation*}
where $\vphi$ is any test function such that $\rmspt T \subset \rmint \{ \vphi = 1 \}$, see \cite[4.1.10]{GMT}.
In particular, if $\la \vphi , T \ra = 0$ then $\partial (\lseg a \rseg \cone T) = T$.
The condition $\la \vphi , T \ra = 0$ for such $\vphi$ is equivalent to $\chi(T) = 0$ in terms of the augmentation map $\chi$, see \ref{P.AM}.
\end{Empty}

\begin{Empty}[Flat norm and flat chains]
\label{P.FNFC}
For $T \in \calD_m(\Rn)$ we define the {\em flat norm}\footnote{Our version differs from $\bF_K$ defined in \cite[4.1.12]{GMT} in that we do not impose restrictions on support.} of $T$ by the formula $\bF(T) = \inf \{ \bM(R) + \bM(S) : T = R + \partial S, R \in \calD_m(\Rn) \text{ and } S \in \calD_{m+1}(\Rn) \} \in [0, \infty]$.
The flat norm was introduced by H. Whitney \cite{WHITNEY}.
It is a norm, by the same argument as in \cite[4.1.12]{GMT} (akin the proof of \ref{3.2} and \ref{5.2}(A) below).
\par 
It is easy to check that $\calD_m(\Rn) \cap \{ T : \bF(T) < \infty \}$ equipped with the norm $\bF$ is a Banach space that contains $\bN_m(\Rn)$ (since $\bF \leq \bM$).
We then define $\bF_m(\Rn)$ as the $\bF$-closure of $\bN_m(\Rn)$ in the former.
Therefore, $\bF_m(\Rn)[\bF]$ is a Banach space.
Its members are called the {\em $m$-dimensional flat chains in $\Rn$}.
\par 
If $T \in \bF_m(\Rn)$ and $m \geq 1$ then $\partial T \in \bF_{m-1}(\Rn)$.
This follows from the obvious inequality $\bF(\partial T) \leq \bF(T)$ and the last sentence of \ref{P.NC}.
\par 
The flat norm is useful in conjunction with the homotopy formula, see \cite[4.1.9]{GMT}.
We now quote two important consequences of the latter.
Below, we abbreviate $\bF_{m,\rmcpt}(\Rn) = \bF_m(\Rn) \cap \{ T : \rmspt T \text{ is compact}\}$.
\begin{enumerate}
\item[(A)] {\it If $\Rn \to \Rp$ is Lipschitz then there exist $f_\# : \bF_{m,\rmcpt}(\Rn) \to \bF_{m,\rmcpt}(\Rp)$ for all $m$ satisfying the following properties.
\begin{enumerate}
\item[(a)] $f_\#$ is linear and $\bF$-continuous. Specifically, $\bF(f_\#T) \leq C(f) \cdot \bF(T)$, where $C(f) = \max \{ (\rmLip f)^m , (\rmLip f)^{m+1} \}$. Moreover, if $f$ is smooth then $f_\# T$ coincides with the same symbol defined in \ref{P.BPF}.
\item[(b)] $\partial \circ f_\# = f_\# \circ \partial$.
\item[(c)] $(g \circ f)_\# = g_\# \circ f_\#$ for all Lipschitz maps $g : \Rp \to \R^q$.
\item[(d)] $\rmspt f_\# T \subset f(\rmspt T)$.
\item[(e)] $\bM(f_\# T) \leq (\rmLip f)^m \bM(T)$.
\end{enumerate}
}
\end{enumerate}
\par 
This can be proved as in \cite[4.1.14]{GMT} by regularization (see \ref{P.SMOOTH}).
We note that all conclusions hold with $\bF_{m,\rmcpt}(\Rn)$ replaced by $\bF_m(\Rn)$ but neither will we use that nor can we quote a reference.\cqfd
\par 
Possibly one of the main properties of flat chains is the following.
\begin{enumerate}
\item[(B)] {\it If $f,g : \Rn \to \Rp$ are Lipschitz, $T \in \bF_{m,\rmcpt}(\Rn)$, and $f|_{\rmspt T} = g|_{\rmspt T}$ then $f_\#T = g_\#T$.}
\end{enumerate}
\par 
See \cite[4.1.15]{GMT} or \cite[26.24]{SIMON} (the latter deals only with the case when $T$ is normal and $f$ is $C^1$ but the way of using the homotopy formula is there).\cqfd
\par 
It is important to realize that the condition $\bF(T) < \infty$ does not imply that $T$ is a flat chain as the following ``hirsute'' example illustrates.
Let $T \in \calD_1(\R^2)$ be defined by $\la \vphi_1 \cdot dx_1 + \vphi_2 \cdot dx_2 , T \ra = \int_0^1 \vphi_2(t,0)d\calL^1(t)$.
Note that $\bF(T) \leq \bM(T) < \infty$.
Consider the projection on the first axis $f : \R^2 \to \R^2 : (x_1,x_2) \mapsto (x_1,0)$.
Here, we accept that $f_\#T$ defined in (A) satisfies $\la  \vphi_1 \cdot dx_1 + \vphi_2 \cdot dx_2 , f_\# T \ra = \la f^\#(\vphi_1 \cdot dx_1 + \vphi_2 \cdot dx_2),T \ra$.
Since $ f^\#(\vphi_1 \cdot dx_1 + \vphi_2 \cdot dx_2) = (\vphi_1 \circ f) \cdot df_1 + (\vphi_2 \circ f) \cdot df_2 = (\vphi_1 \circ f) \cdot dx_1$, we conclude that $f_\# T = 0$.
On the other hand $\rmid_{\R^2\,\#}T = T \neq 0$.
As $f|_{\rmspt T} = \rmid_{\R^2}|_{\rmspt T}$, it follows from (B) that $T \not \in \bF_1(\R^2)$.
\par 
In fact, it follows from (B), (A)(c), and the characterization of $\bF_1(\R)$ (see \eg \cite[4.1.18 last \S\, of p.376]{GMT} but there are easier ways to do this) that if $T \in \bF_1(\R^2)$ and $\rmspt T \subset \R^2 \cap \{ (x_1,0) : 0 \leq x_1 \leq 1 \}$ then $\la \vphi_1 \cdot dx_1 + \vphi_2 \cdot dx_2 , T \ra = \int_0^1 \vphi_1(t,0) \cdot u(t) \, d\calL^1(t)$ for some Lebesgue-summable function $u : [0,1] \to \R$.
The same holds for all $m$ and for the line segment replaced by lipeomorphic images of open subsets of $\Rm$.
\end{Empty}

\begin{Empty}[$T \hel \psi$]
\label{P.RESTR}
If $T \in \calD_m(\Rn)$, $\psi \in \calE^k(\Rn)$, and $0 \leq k \leq m$ then we define $T \hel \psi \in \calD_{m-k}(\Rn)$ by means of the formula $\la \omega , T \hel \psi \ra = \la \psi \wedge \omega , T \ra$ for all $\omega \in \calD^{m-k}(\Rn)$.
Notice that $(T,\psi) \mapsto T \hel \psi$ is bilinear.
We assume throughout this number that $0 \leq k \leq m$.
\begin{enumerate}
\item[(A)] {\it If $T \in \calD_m(\Rn)$ and $\psi \in \calE^{k}(\Rn)$ then $\rmspt T \hel \psi \subset (\rmspt \psi) \cap (\rmspt T)$.}
\end{enumerate}
\par 
This follows easily from the observation that $\rmspt(\psi \wedge \omega) \subset (\rmspt \psi) \cap (\rmspt \omega)$ for all $\omega \in \calD^{k}(\Rn)$.\cqfd
\begin{enumerate}
\item[(B)] {\it If $T \in \calD_m(\Rn)$, $\psi \in \calD^{0}(\Rn)$, and $\rmspt T \subset \rmint \{ \psi = 1 \}$ then $T \hel \psi = T$.}
\end{enumerate}
\par 
We infer from (A) that $\rmspt T \hel (\psi - \ind_{\Rn}) \subset \rmspt (T) \cap \rmspt (\psi - \ind_{\Rn}) = \emptyset$.
Thus, $T \hel (\psi - \ind_{\Rn}) = 0$.\cqfd
\begin{enumerate}
\item[(C)] {\it If $k < m$, $T \in \calD_m(\Rn)$, and $\psi \in \calE^{k}(\Rn)$ then $(-1)^k\partial(T \hel \psi) = (\partial T) \hel \psi - T \hel d\psi$.}
\end{enumerate}
\par 
This is a consequence of the equation $d(\psi \wedge \omega) = (d\psi) \wedge \omega + (-1)^k \psi \wedge d\omega$.\cqfd
\begin{enumerate}
\item[(D)] {\it If $k=0$ or $k=1$, $T \in \calD_m(\Rn)$, $\psi \in \calD^{k}(\Rn)$, and $\bM(T) < \infty$ then $\bM(T \hel \psi) \leq \bM(\psi) \cdot \|T\|\{ \psi \neq 0 \} \leq \bM(\psi) \cdot \bM(T)$.}
\end{enumerate}
\par 
Let $\omega \in \calD^{m-k}(\Rn)$ and $x \in \Rn$.
Since $\psi(x) \in \bigwedge^k\Rn$ and $k=0$ or $k=1$, we have $|\psi(x) \wedge \omega(x)| \leq |\psi(x)| \cdot | \omega(x) |$ (if $k > 1$ then a factor depending on dimensions appears on the right side).
The conclusion now easily follows from the integral representation of $T$, see \ref{P.M}.\cqfd
\begin{enumerate}
\item[(E)] {\it If $T \in \calD_m(\Rn)$, $\psi \in \calD^0(\Rn)$, and $\bN(T) < \infty$ then $T \hel \psi \in \bN_m(\Rn)$.}
\end{enumerate}
\par 
It follows from (A) that $T \hel \psi$ has compact support and it follows from (D) that $\bM(T \hel \psi) < \infty$.
Finally, if $m \geq 1$ then (C) and (D) imply that $\bM [ \partial (\psi \hel \psi) ] < \infty$.\cqfd
\begin{enumerate}
\item[(F)] {\it If $T \in \calD_m(\Rn)$ and $\psi \in \calD^{0}(\Rn)$ then $\bF(T \hel \psi) \leq (\bM(\psi) + \bM(d\psi)) \cdot \bF(T)$.}
\end{enumerate}
\par 
Let $R \in \calD_m(\Rn)$ and $S \in \calD_{m+1}(\Rn)$ be so that $T = R + \partial S$.
Then $T \hel \psi = R' + \partial S'$, where $R' = R \hel \psi + S \hel d\psi$ and $S' = S \hel \psi$, according to (C).
Whence, $\bF(T \hel \psi) \leq \bM(R') + \bM(S') \leq \bM(\psi) \cdot \bM(R) + \bM(d\psi) \cdot \bM(S) + \bM(\psi) \cdot \bM(S) \leq (\bM(\psi) + \bM(d\psi)) \cdot (\bM(R) + \bM(S))$, by (D).
Taking the infimum over all such pairs $(R,S)$ yields the conclusion.\cqfd
\end{Empty}

The following provides the link between our definition of flat norm and that of \cite[4.1.7]{GMT}.
The main point is that if $T$ is a compactly supported flat chain then, in the definition of $\bF(T)$ (recall \ref{P.FNFC}), one can restrict to $R$ and $S$ having compact support as well.
We will refer to this in section \ref{sec.rep}.

\begin{Theorem}
\label{P.THM}
Let $T \in \bF_m(\Rn)$ be such that $\rmspt T \subset \bC^n$ and $\veps > 0$.
Then there exist $R \in \calD_m(\Rn)$ and $S \in \calD_{m+1}(\Rn)$ such that 
\begin{enumerate}
\item[(A)] $T = R + \partial S$;
\item[(B)] $\bM(R) + \bM(S) < \veps + \bF(T)$;
\item[(C)] $(\rmspt R) \cup (\rmspt S) \subset \bC^n$.
\end{enumerate}
\end{Theorem}

\begin{proof}
{\bf (i)} 
{\it The space $\bN_m(\Rn) \cap \{ T : \rmspt T \subset \bC^n \}$ is $\bF$-dense in $\bF_m(\Rn) \cap \{ T : \rmspt T \subset \bC^n\}$.}
\par 
Let $T \in \bF_m(\Rn)$ be supported in $\bC^n$.
By definition of $\bF_m(\Rn)$, there is a sequence $\la T_j \ra_j$ in $\bN_m(\Rn)$ such that $\bF(T-T_j) \to 0$ as $j \to \infty$.
Note that $T = \pi_\#T$, according to \ref{P.FNFC}(B), since $\pi|_{\bC^n} = \rmid_{\Rn}|_{\bC^n}$ and $\rmspt T \subset \bC^n$.
Therefore $\bF(T-\pi_\#T_j) = \bF(\pi_\#(T-T_j)) \leq \bF(T-T_j)$, by \ref{P.FNFC}(A)(a).
Since $\rmspt \pi_\# T_j \subset \bC^n$, by \ref{P.FNFC}(A)(d), and $\bN(\pi_\#T_j) \leq \bN(T_j)$, the proof of {\bf (i)} is complete.
\par 
{\bf (ii)}
Applying Hahn-Banach's theorem as in \cite[4.1.12 pp.367-68]{GMT} (see also \ref{3.2} for a similar argument) we infer that for each $T \in \calD_m(\Rn)$ with $\bF(T) < \infty$ there exist $R \in \calD_m(\Rn)$ and $S \in \calD_{m+1}(\Rn)$ such that $T = R + \partial S$ and $\bM(R) + \bM(S) = \bF(T)$.
\par 
Furthermore, if $T$ is normal then $\bN(R) + \bN(S) < \infty$.
Indeed, $\bM(S) \leq \bF(T) < \infty$ and $\bM(\partial S) \leq \bM(T) + \bM(R) \leq \bM(T) + \bF(T) < \infty$.
Moreover, $\bM(R) \leq \bF(T) < \infty$ and if $m \geq 1$ then $\partial R = \partial T$ so that $\bM(\partial R) = \bM(\partial T) < \infty$.
\par 
{\bf (iii)}
{\it If $T \in \bN_m(\Rn)$, $\rmspt T \subset \bC^n$, and $\hat{\veps} > 0$ then there exist $R \in \bN_m(\Rn)$ and $S \in \bN_{m+1}(\Rn)$ such that $T = R + \partial S$, $\bM(R) + \bM(S) < \hat{\veps} + \bF(T)$, and $(\rmspt R) \cup (\rmspt S) \subset \bC^n$.}
\par 
Given such a $T$ we infer from {\bf (ii)} that there exist $R \in \calD_m(\Rn)$ and $S \in \calD_{m+1}(\Rn)$ such that $T = R + \partial S$ and $\bN(R) + \bN(S) < \infty$.
Letting $\psi \in \calD^0(\Rn)$ be such that $\bC^n \subset \rmint \{ \psi = 1 \}$ recall that $T = T \hel \psi$, by \ref{P.RESTR}(B), and note that $T = R' + \partial S'$ with the same notations as in the proof of \ref{P.RESTR}(F).
Moreover, $R' \in \bN_m(\Rn)$ and $S' \in \bN_{m+1}(\Rn)$, according to \ref{P.RESTR}(E).
We further assume that $\bM(\psi)=1$ so that $\bM(R') + \bM(S')\leq (1 + \rmLip \psi) \cdot \bF(T)$, as in the proof of \ref{P.RESTR}(F).
Since $\rmLip \psi$ can be made arbitrarily small, we can select $\psi$ so that $\bM(R') + \bM(S') < \hat{\veps} + \bF(T)$.
Finally, as in the proof of {\bf (i)}, we note that $T = \pi_\#T = \pi_\#R' + \partial \pi_\#S'$, $\bM(\pi_\#R') + \bM(\pi_\#S') \leq \bM(R') + \bM(S') < \hat{\veps} + \bF(T)$ and $\pi_\#R'$ (resp. $\pi_\#S'$) is an $m$-dimensional (resp. $(m+1)$-dimensional) normal current supported in $\bC^n$.
\par 
{\bf (iv)}
We turn to proving the theorem. 
Let $T \in \bF_m(\Rn)$ be supported in $\bC^n$ and let $\veps > 0$.
It follows from {\bf (i)} that there exists a sequence $\la T_j \ra_j$ of $m$-dimensional normal currents in $\Rn$, all supported in $\bC^n$, such that $T = \sum_{j=1}^\infty T_j$ ($\bF$-convergent series) and $\sum_{j=1}^\infty \bF(T_j) < \frac{\veps}{2} + \bF(T)$.
For each $j$ we apply {\bf (iii)} to $T_j$ and $\hat{\veps} = \frac{\veps}{2^{j+1}}$ and we find $R_j \in \bN_m(\Rn)$ and $S_j \in \bN_{m+1}(\Rn)$ such that $T_j = R_j + \partial S_j$, $\bM(R_j) + \bM(S_j) < \frac{\veps}{2^{j+1}} + \bF(T_j)$, and $(\rmspt R_j) \cup (\rmspt S_j) \subset \bC^n$.
It readily ensues that $\sum_{j=1}^\infty \bM(R_j) < \infty$ and $\sum_{j=1}^\infty \bM(S_j) < \infty$.
Therefore, there exist $R \in \calD_m(\Rn)$ and $S \in \calD_{m+1}(\Rn)$ such that $R = \sum_{j=1}^\infty R_j$ and $S = \sum_{j=1}^\infty S_j$ (both weakly*-convergent series).
It follows that $T = \sum_{j=1}^\infty T_j = \sum_{j=1}^\infty (R_j + \partial S_j) = \sum_{j=1}^\infty R_j + \partial \sum_{j=1}^\infty S_j = R + \partial S$.
Moreover, $\bM(R) + \bM(S) \leq \sum_{j=1}^\infty (\bM(R_j) + \bM(S_j)) \leq \frac{\veps}{2} + \sum_{j=1}^\infty \bF(T_j) < \veps + \bF(T)$.
Finally, it is obvious that $(\rmspt R) \cup (\rmspt S) \subset \bC^n$.
\end{proof}

\begin{Theorem}
\label{P.SUPP}
If $T \in \bF_m(\Rn)$ has compact support, $U \subset \Rn$ is open, and $\rmspt T \subset U$ then there exists a sequence $\la T_j \ra_j$ in $\bN_m(\Rn)$ such that $\bF(T-T_j) \to 0$ as $j \to \infty$ and $\rmspt T_j \subset U$ for all $j$.
\end{Theorem}

\begin{proof}
Choose $\psi \in \calD^0(\Rn)$ such that $\rmspt T \subset \rmint \{ \psi = 1 \}$ and $\rmclos(\rmspt \psi) \subset U$ and choose a sequence $\la T_j \ra_j$ in $\bN_m(\Rn)$ such that $\bF(T-T_j) \to 0$ as $j \to \infty$.
Observe that each $T_j \in \bN_m(\Rn)$ and is supported in $U$, by \ref{P.RESTR}(E,A).
As $T = T \hel \psi$, according to \ref{P.RESTR}(B), we infer from \ref{P.RESTR}(F) that $\bF(T - T_j \hel \psi) = \bF[(T-T_j) \hel \psi] \leq (\bM(\psi)+\rmLip \psi) \cdot \bF(T-T_j) \to 0$ as $j \to \infty$.
\end{proof}

\begin{Empty}[Augmentation map]
\label{P.AM}
The content of this number seems to be new.
We define
\begin{equation*}
\chi : \bN_0(\Rn) \to \R : T \mapsto \mu_T(\Rn),
\end{equation*}
recall \ref{P.M}.
\begin{enumerate}
\item[(A)] {\it For all $T \in \bN_0(\Rn)$ and all $\vphi \in \calD^0(\Rn)$ such that $\rmspt T \subset \rmint \{ \vphi = 1 \}$ we have $\chi(T) = \la \vphi , T \ra$.}
\end{enumerate}
\par 
Simply note that $\mu_T(\Rn) = \int_{\Rn} \ind_{\Rn} d\mu_T = \int_{\Rn} \vec{T} d\|T\| = \int_{\Rn} (\vec{T} \cdot \vphi) d\|T\| = \la \vphi , T \ra$.\cqfd
\begin{enumerate}
\item[(B)] {\it $\chi$ is linear and $|\chi(T)| \leq \bM(T)$ for all $T \in \bN_0(\Rn)$.}
\end{enumerate}
\par 
Both conclusions follow from (A) upon recalling that normal currents have compact support so that given $T \in \bN_0(\Rn)$ there exists $\vphi \in \calD^0(\Rn)$ with $\rmspt T \subset \rmint \{ \vphi = 1 \}$.\cqfd
\begin{enumerate}
\item[(C)] {\it If $S \in \calD_1(\Rn)$ has compact support, $\vphi \in \calD^0(\Rn)$, and $\rmspt S \subset \rmint \{ \vphi =1 \}$ then $\la \vphi , \partial S \ra = 0$. In particular, if $\partial S \in \bN_0(\Rn)$ then $\chi(\partial S) = 0$.}
\end{enumerate}
\par 
Since $(\rmspt d\vphi) \cap (\rmspt S) = \emptyset$, we have $\la \vphi , \partial S \ra = \la d\vphi , S \ra = 0$.
The second part of the conclusion follows from (A).\cqfd
\begin{enumerate}
\item[(D)] {\it For all $T \in \bN_0(\Rn)$ we have $|\chi(T)| \leq \bF(T)$.}
\end{enumerate}
\par 
Let $\veps > 0$.
Reasoning as in the proof of \ref{P.THM}{\bf (iii)} we find $R \in \bN_0(\Rn)$ and $S \in \bN_1(\Rn)$ such that $T = R + \partial S$ and $\bM(R) + \bM(S) < \veps + \bF(T)$.
Thus, $\chi(T) = \chi(R + \partial S) = \chi(R) + \chi(\partial S) = \chi(R)$, by (B) and (C).
Accordingly, $|\chi(T)| = |\chi(R)| \leq \bM(R) < \veps + \bF(T)$, by (B).
Since $\veps$ is arbitrary, the proof is complete.\cqfd
\par 
If follows from (B) and (D) that $\chi$ extends uniquely to a linear $\bF$-continuous map
\begin{equation*}
\chi : \bF_0(\Rn) \to \R
\end{equation*}
which we call the {\em augmentation map}.
\begin{enumerate}
\item[(E)] {\it If $T \in \bF_0(\Rn)$ has compact support and $\vphi \in \calD^0(\Rn)$ is such that $\rmspt T \subset \rmint \{ \vphi = 1 \}$ then $\chi(T) = \la \vphi , T \ra$.}
\end{enumerate}
\par 
Let $U = \rmint \{ \vphi = 1 \}$.
Referring to \ref{P.SUPP} we may choose a sequence $\la T_j \ra_j$ in $\bN_0(\Rn)$ that is $\bF$-convergent to $T$ and so that $\rmspt T_j \subset U$ for all $j$.
Note that $\chi(T_j) = \la \vphi , T_j \ra$, by (A).
Since $| \chi(T) - \chi(T_j) | \leq \bF(T-T_j) \to 0$ and $\la \vphi , T_j \ra \to \la \vphi , T \ra$ as $j \to \infty$, the proof is complete.\cqfd
\begin{enumerate}
\item[(F)] {\it If $S \in \calD_1(\Rn)$ has compact support and $\partial S \in \bF_0(\Rn)$ then $\chi(\partial S)=0$.}
\end{enumerate}
\par 
This follows from (E) and (C).\cqfd
\end{Empty}

\begin{Empty}[Summable currents, smooth currents, and smoothing]
\label{P.SMOOTH}
An $m$-dimensional current $T$ in $\Rn$ is called {\em summable} if and only if there exists a Lebesgue-summable multi-vector field $\eta : \Rn \to \bigwedge_m \Rn$ such that $\la \omega , T \ra = \int_{\Rn} \la \eta , \omega \ra d\calL^n$ for all $\omega \in \calD^m(\Rn)$.
In that case, we write $T = \calL^n \wedge \eta$.
Note that $\bM(\calL^n \wedge \eta) = \int_{\Rn} |\eta|d\calL^n$.
\par 
A current $T \in \calD_m(\Rn)$ is called {\em smooth} if and only if there exists $\eta \in \calE^m(\Rn)$ such that $\la \omega , T \ra = \int_{\Rn} \la \eta , \omega \ra d\calL^n$ for all $\omega \in \calD^m(\Rn)$.
\begin{enumerate}
\item[(A)] {\it If $T$ is smooth and has compact support then $T$ is normal.}
\end{enumerate}
\par 
If $T$ is smooth and has compact support then $\bM(T) = \int_{\Rn} |\eta| d\calL^n < \infty$.
Moreover, if $m \geq 1$ then $\partial (\calL^n \wedge \eta) = \calL^n \wedge (\rmdiv \eta)$, where $\rmdiv \eta = \sum_{i=1}^n (D_i \eta) \hel dx_i \in \calD^{m-1}(\Rn)$, as follows from integration by parts and the formula $d\omega = \sum_{i=1}^n dx_i \wedge D_i \omega$.\cqfd
\par 
Let $\la \phi_\veps \ra_{\veps > 0}$ be an approximation of the identity such that $\phi_\veps(x) = \veps^{-n} \cdot \phi \left( \frac{x}{\veps} \right)$ for some non-negative even test function $\phi$ with $\int_{\Rn} \phi d\calL^n = 1$ and $\rmspt \phi \subset \bB(0,1)$.
Thus, $\rmspt \phi_\veps \subset \bB(0,\veps)$.
If $T \in \calD_m(\Rn)$ then we define $\phi_\veps * T$ be means of the formula $\la \omega, \phi_\veps * T \ra = \la \phi_\veps * \omega , T \ra$ for all $\omega \in \calD^m(\Rn)$.
One checks that $\phi_\veps * T \in \calD_m(\Rn)$.
The following hold.
\begin{enumerate}
\item[(B)] {\it $\phi_\veps * T$ is smooth.}
\item[(C)] $\rmspt \phi_\veps * T \subset \Rn \cap \{ x : \rmdist(x , \rmspt T) \leq \veps \}$.
\item[(D)] $\bM(\phi_\veps * T) \leq \bM(T)$.
\end{enumerate}
\par 
Conclusions (B) and (C) are standard whereas (D) immediately follows from the observation that $\bM(\psi_\veps * \omega) \leq \bM(\omega)$ for all $\omega$.\cqfd
\begin{enumerate}
\item[(E)] $\bF(T - \phi_\veps * T) \leq \veps \cdot \bN(T)$.
\end{enumerate}
\par 
The case $m=0$ is a trivial consequence of the definitions and the observation that $\bF(T) = \sup \{ \la \omega , T \ra : \omega \in \calD^m(\Rn) \text{ and } \max \{ \bM(\omega) , \bM(d\omega) \} \leq 1 \}$.
The latter can be proved as in \cite[4.1.12]{GMT}.
However, we will also apply (E) to the case $m=1$ whose proof is less trivial. 
It relies on a homotopy formula for the so-called smoothing deformation chain, see the beginning of \cite[4.1.18]{GMT}.\cqfd
\begin{enumerate}
\item[(F)] {\it If $T$ is summable and has compact support then $T \in \bF_m(\Rn)$.}
\end{enumerate}
\par 
Let $\eta : \Rn \to \bigwedge_m \Rn$ be Lebesgue-summable and such that $T = \calL^n \wedge \eta$.
Note that $\phi_\veps * T = \calL^n \wedge (\phi_\veps * \eta)$, whence, $\bM(\phi_\veps * T - T) = \int_{\Rn} |\eta - \phi_\veps * \eta\| d\calL^n$.
Therefore, $\bF(T - \phi_\veps * T) \leq \bM(T - \phi_\veps * T) \to 0$ as $\veps \to 0^+$.
It remains to recall (A) that each $\phi_\veps * T$ is normal.\cqfd
\end{Empty}

\section{The norm $\bG$}

\begin{Empty}[Definition of $\bG$]
\label{3.1}
With a distribution $T \in \calD_0(\Rn)$ we associate 
\begin{equation*}
\bG(T) = \inf \{ \bM(S) : S \in \calD_1(\Rn) \text{ such that } T = \partial S \} \in [0,\infty].
\end{equation*}
In particular, $\bG(T) = \infty$ in case $T$ is not the boundary of any 1-dimensional current.
Note that $\bG(T_1 + T_2) \leq \bG(T_1) + \bG(T_2)$ for all distributions $T_1,T_2$ and that $\bG(\lambda \cdot T) = |\lambda| \cdot \bG(T)$ for all distribution $T$ and real number $\lambda \neq 0$.
The following two straightforward observations will prove useful.
\begin{enumerate}
\item[(A)] {\it For every $T \in \calD_0(\Rn)$ we have $\bF(T) \leq \bG(T)$.}
\end{enumerate}
\par 
If $T = \partial S$ then $\bF(T) \leq \bM(S)$.
Taking the infimum on the right side over all such $S$ yields the conclusion.\cqfd
\begin{enumerate}
\item[(B)] {\it For every $S \in \calD_1(\Rn)$ we have $\bG(\partial S) \leq \bF(S)$.}
\end{enumerate}
\par 
If $S = A + \partial B$ then $\partial S = \partial A$, whence, $\bG(\partial S) \leq \bM(A) \leq \bM(A) + \bM(B)$.
Taking the infimum on the right side over all such couples $(A,B)$ yields the conclusion.\cqfd
\end{Empty}

The following adapts part of \cite[4.1.12]{GMT} to the present situation.

\begin{Theorem}
\label{3.2}
Let $T \in \calD_0(\Rn)$.
The following hold.
\begin{enumerate}
\item[(A)] $\bG(T) = \sup \{ \la \vphi , T \ra : \vphi \in \calD^0(\Rn) \text{ and } \bM(d\vphi) \leq 1 \}$.
\item[(B)] Assume that $\bG(T) < \infty$. Then there exists $S \in \calD_{1}(\Rn)$ such that $T = \partial S$ and $\bG(T) = \bM(S)$. 
\end{enumerate}
\end{Theorem}

\begin{proof}
{\bf (i)}
For the purpose of this proof, we abbreviate
\begin{equation*}
\hat{\bG}(T) = \sup \{ \la \vphi , T \ra : \vphi \in \calD^0(\Rn) \text{ and } \bM(d\vphi) \leq 1 \} \in [0,\infty].
\end{equation*}
\par 
{\bf (ii)}
{\it Here, we show that $\hat{\bG}(T) \leq \bG(T)$.}
\par 
Since this readily holds when $\bG(T) = \infty$, we henceforth assume that $\bG(T) < \infty$.
Let $\vphi \in \calD^0(\Rn)$ and let $S \in \calD_{1}(\Rn)$ be such that $\bM(S) < \infty$ and $T = \partial S$ (note that such $S$ exists, as we assume that $\bG(T) < \infty$).
We have
\begin{equation*}
\la \vphi , T \ra = \la d\vphi , S \ra \leq \bM(d\vphi) \cdot  \bM(S).
\end{equation*}
By the arbitrariness of $\vphi$ and $S$, we conclude that $\hat{\bG}(T) \leq \bG(T)$.
\par 
{\bf (iii)}
{\it Here, we show that if $\hat{\bG}(T) < \infty$ then there exists $S \in \calD_1(\Rn)$ such that $T = \partial S$ and $\bM(S) \leq \hat{\bG}(T)$. In particular, $\bG(T) \leq \hat{\bG}(T)$.}
\par 
We consider the vector space $\calD^{1}(\Rn)$ equipped with the norm $\bM$ as well as the linear map $d : \calD^0(\Rn) \to \calD^{1}(\Rn) : \vphi \mapsto d\vphi$.
We define a function $L$ on $\rmim d$ by the formula $L(d\vphi) = \la \vphi , T \ra$.
Observe that $L$ is unambiguously defined, since $\hat{\bG}(T) < \infty$, and linear.
Notice also that $|L| \leq \hat{\bG}(T) \bM|_{\rmim d}$.
According to the Hahn-Banach theorem, $L$ admits a linear extension $S$ to $\calD^{1}(\Rn)$ satisfying $|S| \leq \hat{\bG}(T) \bM$.
Therefore, $\bM(S) \leq \hat{\bG}(T)$.
In particular, $\bM(S) < \infty$, which implies that $S$ is a $1$-dimensional current in $\Rn$.
Note also that $\partial S = T$.
\par 
{\bf (iv)}
We now establish conclusion (A).
Firstly, $\hat{\bG}(T) \leq \bG(T)$, by {\bf (ii)}.
Secondly, $\bG(T) \leq \hat{\bG}(T)$, regardless of whether $\hat{\bG}(T) < \infty$ or not, by {\bf (iii)}.
\par 
{\bf (v)}
It remains to prove conclusion (B).
If $\bG(T) < \infty$ then $\hat{\bG}(T) < \infty$, according to {\bf (ii)}, whence, there exists $S \in \calD_1(\Rn)$ such that $T = \partial S$ and $\bG(T) \leq \bM(S) \leq \hat{\bG}(T)$, by {\bf (iii)}.
Since $\bG(T) = \hat{\bG}(T)$, (B) holds.
\end{proof}

The following may be interpreted as a linear isoperimetric inequality for 1-dimensional currents.

\begin{Theorem}
\label{3.3}
Let $S \in \calD_1(\Rn)$ have compact support, $T = \partial S$, and $a \in \Rn$.
Then
\begin{equation*}
\bG(T) \leq \max \left\{ |x-a| : x \in \rmspt S \right\} \cdot \bM(T).
\end{equation*}
\end{Theorem}

\begin{proof}
{\bf (i)}
If $\psi \in \calD^0(\Rn)$ and $\rmspt S \subset \rmint \{ \psi = 1 \}$ then $\la \psi, T \ra = \la \psi , \partial S \ra = \la d\psi , S \ra = 0$.
\par 
{\bf (ii)}
With $\vphi \in \calD^0(\Rn)$ we associate $\vphi_a \in \calE^0(\Rn)$ defined by the formula $\vphi_a(x) = \vphi(x) - \vphi(a)$.
Note that $|\vphi_a(x)| = |\vphi(x) - \vphi(a)| \leq (\rmLip \vphi) \cdot |x-a| = \bM(d\vphi) \cdot | x-a | \leq \bM(d\vphi) \cdot \gamma(a,S)$, whenever $x \in \rmspt S$, where $\gamma(a,S) = \max \{ |x-a| : x \in \rmspt S \}$.
\par 
{\bf (iii)}
Let $\psi$ and $\vphi$ be, respectively, as in {\bf (i)} and {\bf (ii)}.
As $(\rmspt \vphi \cdot (\ind_{\Rn} - \psi)) \cap (\rmspt T) = \emptyset$, we have $\la \vphi , T \ra = \la \vphi \cdot \psi , T \ra$ and, by {\bf (i)}, $\la \vphi \cdot \psi , T \ra = \la \vphi \cdot \psi - \vphi(a) \cdot \psi , T \ra = \la \vphi_a \cdot \psi , T \ra \leq \bM(\vphi_a \cdot \psi) \cdot \bM(T)$.
Observe that the infimum of $\bM(\vphi_a \cdot \psi)$ taken over all $\psi$ as in {\bf (i)} equals $\max \{ |\vphi_a(x)| : x \in \rmspt S \}$.
Therefore, $\la \vphi,T\ra \leq \bM(d\vphi) \cdot \gamma(a,S) \cdot \bM(T)$, according to {\bf (ii)}.
Since $\vphi$ is arbitrary, the proof is complete. 
\end{proof}

\begin{Remark}
Note that under the assumptions of \ref{3.3} both $\bG(T)$ and $\bM(T)$ may be equal to $\infty$.
\end{Remark}

\section{Representation of $0$-dimensional flat cycles}
\label{sec.rep}

\begin{Empty}[Definition of $\calG(\bC^n)$]
\label{4.0}
%
%
Recalling the definion of the augmentation map $\chi$ in \ref{P.AM}, we define
\begin{equation*}
\calG(\bC^n) = \bF_0(\Rn) \cap \{ T : \rmspt T \subset \bC^n \text{ and } \chi(T) = 0 \}. 
\end{equation*}
This is a vector subspace of $\calD_0(\Rn)$, the space of distributions in $\Rn$.
We call {\em $0$-dimensional flat cycles in $\bC^n$} the members of $\calG(\bC^n)$.
\end{Empty}

\begin{Theorem}
\label{4.3}
Let $\veps > 0$ and $T \in \calG(\bC^n)$.
There there exists $S \in \bN_1(\Rn)$ such that 
\begin{enumerate}
\item[(A)] $\bG(T - \partial S) < \veps$;
\item[(B)] $\rmspt S \subset \bC^n$.
\end{enumerate}
In particular, $\bG(T) < \infty$.
\end{Theorem}

\begin{proof}
{\bf (i)}
There exists $T_0 \in \bN_0(\Rn)$ such that $\bF(T-T_0) < \frac{\veps}{2}$ and $\rmspt T_0 \subset \bC^n$.
Indeed, the existence of $T_0$ satisfying these properties except perhaps for $\rmspt T_0 \subset \bC^n$ follows from the definition of $\bF_0(\Rn)$, recall \ref{P.FNFC}.
Replacing $T_0$ with $\pi_\#T_0$ if necessary completes the proof of the claim.
\par 
{ \bf (ii)} 
Recalling \ref{P.THM}, we infer that there exist $R_0 \in \calD_0(\Rn)$ and $S_0 \in \calD_1(\Rn)$ such that $T - T_0 = R_0 + \partial S_0$, $\bM(R_0) + \bM(S_0) < \frac{\veps}{2} + \bF(T-T_0) < \veps$ (by {\bf (i)}), and $(\rmspt R_0) \cup (\rmspt S_0) \subset \bC^n$.
\par 
{\bf (iii)}
Note that $T_0 + R_0 \in \bN_0(\Rn) \subset \bF_0(\Rn)$, whence, $\partial S_0 = T - (T_0+R_0) \in \bF_0(\Rn)$.
Therefore, $\chi(\partial S_0)=0$, by \ref{P.AM}(F), and in turn $\chi(T_0+R_0) = \chi(T) + \chi(\partial S_0) = 0$.
\par 
{\bf (iv)}
We define $S = \lseg 0 \rseg \cone (T_0+R_0) \in \bN_1(\Rn)$ so that $\partial S = R_0+T_0$, see \ref{P.CONE}.
Moreover, $\rmspt S \subset \bC^n$, since $\rmspt (T_0 + R_0) \subset \bC^n$ and $\bC^n$ is star-shaped around the origin.
Notice that $T = T_0+R_0+\partial S_0 = \partial S + \partial S_0$.
Therefore, 
\begin{equation*}
\begin{aligned}
\bG(T - \partial S) & = \bG(\partial S_0) && \\
& \leq \bM(S_0) && \text{(by definition of $\bG$)}\\
& < \veps && \text{(by {\bf (ii)}).}
\end{aligned}
\end{equation*}
This completes the proof of conclusions (A) and (B).
Finally,
\begin{equation*}
\bG(T) \leq \bG(T - \partial S) + \bG(\partial S) < \veps + \bM(S) < \infty.
\end{equation*}
\end{proof}

\begin{Theorem}
\label{4.1}
Let $\veps > 0$ and $T \in \calG(\bC^n)$.
Then there exists a \underline{smooth} $1$-dimensional normal current $S \in \bN_1(\Rn)$ such that
\begin{enumerate}
\item[(A)] $\bG(T - \partial S) < \veps$;
\item[(B)] $\bM(S) < \veps + \bG(T)$;
\item[(C)] $\rmspt S \subset \Rn \cap \{ x : \rmdist(x,\bC^n) \leq \veps \}$. 
\end{enumerate}
\end{Theorem}

\begin{proof}
{\bf (i)}
According to \ref{4.3}, there exists $T_0 \in \bN_0(\Rn)$ such that $\bG(T_0) < \infty$, $\chi(T_0) = 0$, $\bG(T-T_0) < \frac{\veps}{3}$, and $\rmspt T_0 \subset \bC^n$.
\par 
{\bf (ii)}
Since $\bG(T_0) < \infty$, there exists $S_0 \in \calD_1(\Rn)$ such that $T_0 = \partial S_0$ and $\bM(S_0) < \frac{\veps}{2} + \bG(T_0)$.
Notice that $\bN(S_0) = \bM(S_0) + \bM(\partial S_0) < \frac{\veps}{2} + \bG(T_0) + \bM(T_0) < \infty$.
We next remedy the fact that $S_0$ may not have compact support.
\par 
{\bf (iii)}
{\it Here, we show that there exists $S_1 \in \bN_1(\Rn)$ such that $\bG(T_0 - \partial S_1) < \frac{\veps}{3}$ and $\bM(S_1) < \frac{\veps}{2} + \bG(T_0)$.}
\par 
Let $1 < r < r + s$ be positive real numbers to be determined momentarily and choose a test function $\psi$ such that $\ind_{\bB(0,r)} \leq \psi \leq \ind_{\bB(0,r+s)}$ and $\rmLip \psi \leq \frac{2}{s}$.
Define $S_1 = S_0 \hel \psi$ and note that $S_1 \in \bN_1(\Rn)$, by \ref{P.RESTR}(E).
Furthermore, $\bM(S_1) \leq \bM(S_0) < \frac{\veps}{2} + \bG(T_0)$, by {\bf (ii)}.
We now turn to showing that $\bG(T_0 - \partial S_1) < \frac{\veps}{3}$ for an appropriate choice of $r$ and $s$.
First, note that $\partial S_1 = (\partial S_0) \hel \psi - S_0 \hel d\psi = T_0 - S_0 \hel d\psi$, in view of \ref{P.RESTR}(C,B).
Accordingly, $\chi(S_0 \hel d\psi) = \chi(T_0 - \partial S_1) = \chi(T_0) - \chi(\partial S_1)=0$, by \ref{P.AM}(C).
Next, define $R = \lseg 0 \rseg \cone (S \hel d\psi) \in \bN_1(\Rn)$ and recall that $\partial R = S \hel d\psi$.
Thus, $T_0 - \partial S_1 = \partial R$ and it follows that $\bG(T_0 - \partial S_1) \leq \bM(R)$.
Now,
\begin{equation*}
\begin{aligned}
\bM(R) & \leq 2 \sqrt{n} \cdot (r+s) \cdot \bM(S_0 \hel d\psi) && \text{(by \ref{P.CONE})} \\
& \leq 2 \sqrt{n} \cdot (r+s) \cdot \bM(d\psi) \cdot \|S_0\| \{ d\psi \neq 0 \} && \text{(by \ref{P.RESTR}(D))} \\
& \leq 2 \sqrt{n} \cdot \left( \frac{r+s}{s}\right) \cdot \|S_0\|(\Rn \setminus \bB(0,r)) .
\end{aligned}
\end{equation*}
One can readily choose $r$ first and $s$ next, both large enough, for $\bM(R) < \frac{\veps}{3}$.
\par 
{\bf (iv)}
{\it Here, we show that there exists $S_2 \in \bN_1(\Rn)$ such that $\bG(T_0 - \partial S_2) < \frac{\veps}{3}$, $\bM(S_2) < \frac{\veps}{2} + \bG(T_0)$, and $\rmspt S_2 \subset \bC^n$.}
\par 
It suffices to put $S_2 = \pi_\# S_1$ which is well-defined, since $S_1$ has compact support.
Thus, $\rmspt S_2 \subset \bC^n$ and the other conclusions follow from {\bf (iii)} and the equation $T_0 = \pi_\#T_0$, recall \ref{P.FNFC}(B).
\par 
{\bf (v)}
We choose $\hat{\veps} > 0$ small enough for $\hat{\veps} \cdot \bN(S_2) < \veps$ and we define $S = \phi_{\hat{\veps}} * S$.
Thus, $S \in \bN_1(\Rn)$ is a smooth current, $\bM(S) \leq \bM(S_2) < \frac{\veps}{2} + \bG(T_0)$, $\rmspt S \subset \Rn \cap \{ x : \rmdist(x,\rmspt S) \leq \veps\}$, and $\bF(S-S_2) \leq \hat{\veps} \cdot \bN(S_2) < \veps$, all by \ref{P.SMOOTH}.
We establish conclusion (A):
\begin{equation*}
\begin{aligned}
\bG(T - \partial S) & \leq \bG(T - T_0) + \bG(T_0 - \partial S_2)  && \\
& \qquad \qquad \qquad \qquad + \bG[\partial (S_2 - S)] && \\
& < \frac{\veps}{3} + \frac{\veps}{3} + \bF(S_2 - S) && \text{(by {\bf (i)}, {\bf (iv)}, and \ref{3.1}(B))} \\
& < \veps . &&
\end{aligned}
\end{equation*}
Regarding conclusion (B), we observe that
\begin{equation*}
\bM(S) < \frac{\veps}{2} + \bG(T_0) \leq \frac{\veps}{2} + \bG(T_0-T) + \bG(T) < \veps + \bG(T),
\end{equation*}
in view of {\bf (i)}.
\end{proof}

In the following statement, the difference with \ref{4.3} is that $\bM(S)$ is now  guaranteed to be close to $\bG(T)$ and the difference with \ref{4.1} is that $S$ is not requested to be smooth but is requested to be supported in $\bC^n$.

\begin{Theorem}
\label{4.4}
Let $\veps > 0$ and $T \in \calG(\bC^n)$.
Then there exists a $1$-dimensional normal current $S \in \bN_1(\Rn)$ such that
\begin{enumerate}
\item[(A)] $\bG(T - \partial S) < \veps$;
\item[(B)] $\bM(S) < \veps + \bG(T)$;
\item[(C)] $\rmspt S \subset \bC^n$. 
\end{enumerate}
\end{Theorem}

\begin{proof}
This has been established in the course of proving \ref{4.1}, namely in step {\bf (iv)} of the proof.
\end{proof}

The following notations will be useful in the proof of the next result.
Given $r > 0$ we let $\bmu_r$ be the proper smooth map $\Rn \to \Rn : x \mapsto r \cdot x$.
Observe that $\bM(\bmu_{r \,\#} T) = r^m \bM(T)$ and $\rmspt \bmu_{r \,\#} T = \bmu_r(\rmspt T)$ for all $T \in \calD_m(\Rn)$.
Furthermore, abbreviate $\bB_\infty(0,r) = \Rn \cap \{ x : |x|_\infty \leq r\}$.
Notice that $\Rn \cap \{ x : \rmdist(x,\bC^n) \leq \veps \} \subset \bB_\infty(0,1+\veps) \subset \Rn \cap \{x : \rmdist(x,\bC^n) \leq \sqrt{n} \cdot \veps\}$.

\begin{Theorem}
\label{4.2}
Let $\veps > 0$ and $T \in \calG(\bC^n)$.
Then there exists a Borel-measurable vector-field $\eta : \Rn \to \bigwedge_1\Rn$ such that:
\begin{enumerate}
\item[(A)] $T = \partial ( \calL^n \wedge \eta)$;
\item[(B)] $\bG(T) \leq \int_{\Rn} |\eta| d\calL^n < \veps + \bG(T)$;
\item[(C)] $\rmspt (\calL^n \wedge \eta) \subset \Rn \cap \{ x : \rmdist(x,\bC^n) \leq \veps \}$.
\end{enumerate}
\end{Theorem}

\begin{proof}
{\bf (i)}
Let $\la \veps_j \ra_{j \geq 1}$ be a null decreasing sequence of positive real numbers whose values will be determined momentarily and define an increasing sequence $\la r_j \ra_{j \geq 0}$ by the requirements that $r_0=1$ and $r_{j} = r_{j-1} + \veps_j$ for all $j \geq 1$.
\par 
{\bf (ii)}
Here, we shall define by induction on $j=1,2,\cdots$ two sequences of currents $\la T_j \ra_{j \geq 0}$ and $\la S_j \ra_{j \geq 1}$ satisfying the following conditions for all $j \geq 1$.
\begin{enumerate}
\item[(1)] $T_{j-1} \in \bF_0(\Rn)$, $\chi(T_{j-1}) = 0$, and $\rmspt T_{j-1} \subset \bB_\infty(0,r_{j-1})$;
\item[(2)] $S_j \in \bN_1(\Rn)$ and is smooth;
\item[(3)] $\bG(T_{j-1} - \partial S_j) < \veps_j$;
\item[(4)] $\bM(S_j) < \veps_j + \bG(T_{j-1})$;
\item[(5)] $\rmspt S_j \subset \bB_\infty(0,r_j)$;
\item[(6)] $T_j = T_{j-1} - \partial S_j$.
\end{enumerate}
\par 
To start the induction, we put $T_0=T$ so that condition (1) is satisfied for $j=1$.
Letting $j \geq 1$ and assuming that (1) holds, we define $\hat{T}_{j-1} = \bmu_{r_{j-1}^{-1} \, \#} T_{j-1}$ and $\hat{\veps}_j = r_{j-1}^{-1}\veps_j$ and we note that \ref{4.1} applies to $\hat{T}_{j-1}$ and $\hat{\veps}_j$ and yields a smooth current $\hat{S}_j \in \bN_1(\Rn)$ such that $\bG ( \hat{T}_{j-1} - \partial \hat{S}_j ) < \hat{\veps}_j$, $\bM ( \hat{S}_j ) < \hat{\veps}_j + \bG ( \hat{T}_{j-1} )$ and $\rmspt \hat{S}_j \subset \bB_\infty \left( 0 , 1 + \hat{\veps}_j \right)$.
Letting $S_j = \bmu_{r_{j-1} \, \#} S_j$, we note that $S_j$ satisfies (2), (3), and (4), since $r_{j-1} \hat{\veps}_j = \veps_j$.
Moreover, $\rmspt S_j \subset \bB_\infty \left( 0, r_{j-1} \left( 1 + \hat{\veps}_j\right)\right)$ and, as $r_{j-1} \left( 1 + \hat{\veps}_j\right) = r_j$, we conclude that $S_j$ satisfies (5).
In turn, we define $T_j$ according to (6) and we observe that $T_j \in \bF_0(\Rn)$, $\chi(T_j) = 0$, by (1) and \ref{P.AM}(C), and $\rmspt T_j \subset (\rmspt T_{j-1}) \cup (\rmspt S_j) \subset \bB_\infty(0,r_j)$, since $\la r_j \ra_j$ is increasing.
\par 
{\bf (iii)}
We observe that
\begin{equation}
\label{eq.1}
T_j = T - \sum_{i=1}^j \partial S_j
\end{equation}
and 
\begin{equation}
\label{eq.2}
\sum_{i=1}^j \bM(S_i) < \veps_j + 2 \cdot \sum_{i=1}^{j-1} \veps_i + \bG(T)
\end{equation}
for all $j \geq 1$.
Both statements are obtained by induction on $j$.
A proof of the former is based on (6) and a proof of the latter is based on (4), (6), and (3).
\par 
{\bf (iv)}
Let $\eta_j : \Rn \to \bigwedge_{1} \Rn$ be compactly supported $C^\infty$ vector-fields such that $S_j = \calL^n \wedge \eta_j$ for all $j \geq 1$.
Recall that $\bM(S_j) = \int_{\Rn} |\eta_j| d\calL^n$.
Choose $\veps_i = \frac{\veps}{\sqrt{n} \cdot 2^{i+2}}$.
It follows from \eqref{eq.2} that the series $\sum_{i=1}^\infty S_i$ is $\bM$-convergent.
Therefore, its sum $S \in \calD_1(\Rn)$ is represented as $S = \calL^n \wedge \eta$ for some Borel-measurable vector-field $\eta : \Rn \to \bigwedge_1 \Rn$ satisfying
\begin{equation*}
\int_{\Rn} |\eta| d\calL^n = \bM(S) \leq \sum_{i=1}^\infty \bM(S_i)  < \veps + \bG(T),
\end{equation*}
by \eqref{eq.2}.
Thus, the second inequality in conclusion (B) holds.
Conclusion (C) holds as well, by (5) and $r_j \leq 1 + \frac{\veps}{\sqrt{n}}$ for all $j \geq 1$.
Finally, we observe that $\bG(T_j) \to 0$ as $j \to \infty$, according to conditions (6) and (3).
Considering $\bG$-convergence in \eqref{eq.1} as $j \to \infty$, we conclude that $T = \partial S$.
\end{proof}

\section{The Banach space $\calG(\bC^n)$}

\begin{Empty}[Abbreviations and definition]
\label{5.1}
We abbreviate
\begin{equation*}
\bM\bF_1(\Rn) = \bF_1(\Rn) \cap \{ S : \bM(S) < \infty \}
\end{equation*}
and 
\begin{equation*}
\bM\bF_{1,\bC^n}(\Rn) = \bM\bF_1(\Rn) \cap \{ S : \rmspt S \subset \bC^n \}
\end{equation*}
and we define 
\begin{equation*}
\bG^\times(T) = \inf \{ \bM(S) : T = \partial S \text{ and } S \in \bM\bF_{1,\bC^n}(\Rn) \} \in [0,\infty]
\end{equation*}
whenever $T \in \calD_0(\Rn)$.
Clearly, both $\bM\bF_1(\Rn)$ and $\bM\bF_{1,\bC^n}(\Rn)$ are vector spaces.
Since $\bF(S) \leq \bM(S)$ for all $S \in \calD_1(\Rn)$, one easily checks that both $\bM\bF_1(\Rn)[\bM]$ and $\bM\bF_{1,\bC^n}(\Rn)[\bM]$ are Banach spaces.
\end{Empty}

\begin{Theorem}
\label{5.2}
The following hold.
\begin{enumerate}
\item[(A)] $\calG(\bC^n)[\bG]$ is a Banach space.
\item[(B)] This is a short exact sequence:
\begin{equation*}
\begin{CD}
\ker \partial @>{\iota}>> \bM\bF_{1,\bC^n}(\Rn)[\bM] @>{\partial}>> \calG(\bC^n)[\bG] @>{\chi}>> \R,
\end{CD}
\end{equation*}
where $\iota$ is the inclusion, \ie $\partial$ is surjective.
\item[(C)] For every $T \in \calG(\bC^n)$ we have $\bG(T) = \bG^\times(T)$.
\end{enumerate}
\end{Theorem}

\begin{proof}[Proof of (A)]
{\bf (i)} 
Recall \ref{4.3} that $\bG(T) < \infty$ whenever $T \in \calG(\bC^n)$.
It is readily observed that $\bG$ is a semi-norm.
If $\bG(T) = 0$ then $T = 0$, in view of \ref{3.2}(A).
Thus, $\bG$ is a norm.
\par 
{\bf (ii)}
Here again, we use the expression of $\bG(T)$ given in \ref{3.2}(A).
Let $\la T_j \ra_j$ be a $\bG$-Cauchy sequence in $\calG(\bC^n)$.
As $\la \la \vphi , T_j \ra \ra_j$ is a Cauchy sequence in $\R$ for every $\vphi \in \calD^0(\Rn)$, the sequence $\la T_j \ra_j$ weakly*-converges to a distribution $T$.
In particular, $\rmspt T \subset \bC^n$.
One routinely checks that $\bG(T-T_j) \to 0$ as $j \to \infty$.
According to \ref{3.1}(A), $\bF(T-T_j) \to 0$ as $j \to \infty$.
Thus, $T \in \bF_0(\Rn)$.
Finally, $\chi(T) = 0$, by the paragraph following \ref{P.AM}(D).
\end{proof}

\begin{proof}[Proof of (B) and (C)]
Let $T \in \calG(\bC^n)$.
It is most obvious that $\bG(T) \leq \bG^\times(T)$.
Let $\veps > 0$ and let $\eta$ be associated with $T$ and $\veps$ in \ref{4.4}.
Recall \ref{P.SMOOTH}(F) that $\calL^n \wedge \eta \in \bF_1(\Rn)$.
As $\rmspt (\calL^n \wedge \eta)$ is compact, we may define $S = \pi_\# (\calL^n \wedge \eta) \in \bF_1(\Rn)$.
Notice that $\rmspt S \subset \bC^n$ and $\bM(S) \leq \bM(\calL^n \wedge \eta) < \veps + \bG(T)$.
In particular, $\bM(S) < \infty$, according to \ref{4.3}.
Thus, $S \in \bM\bF_{1,\bC^n}(\Rn)$.
Since $T = \pi_\#T = \partial \pi_\# ( \calL^n \wedge \eta) = \partial S$,  conclusion (B) follows from the arbitrariness of $T$.
Furthermore, $\bG^\times(T) \leq \bM(S) < \veps + \bG(T)$ so that conclusion (C) follows from the arbitrariness of $\veps$.
\end{proof}

\section{A dense subspace of $\calG(\bC^n)$}

\begin{Empty}[0-dimensional polyhedral currents]
\label{6.1}
We say that $T \in \calD_0(\Rn)$ is a {\em 0-dimen\-sional polyhedral current in $\Rn$} if there are a finite set $F \subset \Rn$ and $\theta : F \to \R$ such that $T = \sum_{x \in F} \theta(x) \cdot \lseg x \rseg$.
In that case, $\rmspt T \subset F$, $T \in \bN_0(\Rn)$, and $\chi(T) = \sum_{x \in F} \theta(x)$.
We let $\bP_0(\Rn)$ be the vector space of 0-dimensional polyhedral currents in $\Rn$.
Notice that $\bP_0(\Rn) \subset \bN_0(\Rn) \subset \bF_0(\Rn)$.
\end{Empty}

\begin{Theorem}
\label{6.2}
The vector subspace
\begin{equation*}
\bP_0(\Rn) \cap \{ P : \rmspt P \subset \bC^n \text{ and } \chi(P) = 0 \}
\end{equation*}
is $\bG$-dense in $\calG(\bC^n)$.
\end{Theorem}

\begin{proof}
{\bf (i)}
Let $\veps > 0$ and $T \in \calG(\bC^n)$.
We ought to prove the existence of $P \in \bP_0(\Rn)$ such that $\rmspt P \subset \bC^n$, $\chi(P) = 0$, and $\bG(T-P) < \veps$.
We choose $T_0$, $R_0$, and $S_0$ as in the proof of \ref{4.3} {\bf (i)} and {\bf (ii)}.
We abbreviate $R = T_0 + R_0$ and we observe that $R \in \bN_0(\Rn)$, $\rmspt R \subset \bC^n$, and $\chi(R) = 0$, by \ref{4.3} {\bf (iii)}.
We also record the inequality $\bG(T - R) = \bG(T - (T_0 + R_0)) = \bG(\partial S_0) \leq \bM(S_0) < \veps$, by the definition of $\bG$ and \ref{4.3} {\bf (ii)}.
\par 
{\bf (ii)}
We choose $\hat{\veps} > 0$ such that $\hat{\veps} \cdot \bM(R) < \veps$ and a finite partition $\calQ$ of $\bC^n$ whose members are Borel-measurable and of diameter less than $\hat{\veps}$.
Furthermore, for each $Q \in \calQ$ we choose $x_Q \in Q$.
Define 
\begin{equation*}
P = \sum_{Q \in \calQ} \chi(R \hel Q) \lseg x_Q \rseg \in \bP_0(\Rn),
\end{equation*}
where $R \hel Q$ are the 0-dimensional normal currents characterized by the equations $\la \vphi , R \hel Q \ra = \int_Q \la \vec{R} , \vphi \ra d\|R\|$.
Note that $\rmspt P \subset \bC^n$ and $\chi(P) = \sum_{Q \in \calQ} \chi(R \hel Q) = \chi(R) = 0$.
\par 
Let $\vphi \in \calD^0(\Rn)$ be such that $\bM(d\vphi) \leq 1$.
We have,
\begin{multline*}
\la \vphi , R - P \ra  = \int_{\Rn} \la \vphi , \vec{R} \ra d\|R\| - \sum_{Q \in \calQ} \chi(R \hel Q) \vphi(x_Q) \\ = \sum_{Q \in \calQ} \int_{\Rn} \la \vphi(x) - \vphi(x_Q) , \vec{R} \ra d\|R\| 
\leq \sum_{Q \in \calQ} (\rmLip \vphi)\cdot (\rmdiam Q) \cdot \|R\|(Q) \\ \leq \hat{\veps} \cdot \bM(R) < \veps.
\end{multline*}
Since $\vphi$ is arbitrary, $\bG(R-P) \leq \veps$, according to \ref{3.2}(A).
Therefore, $\bG(T - P) < 2 \cdot \veps$, by {\bf (i)}.
\end{proof}

\section{$\calF(\bC^n)$ and $\calG(\bC^n)$ are isometrically isomorphic}

\begin{Theorem}
\label{7.1}
The Banach spaces $\calF(\bC^n)[\|\cdot\|_L^*]$ and $\calG(\bC^n)[\bG]$ are isometrically isomorphic.
\end{Theorem}

\begin{proof}
{\bf (i)} 
Recall \ref{P.LF}.
Each member $\alpha$ of the vector subspace $E \subset \rmLip_o(\bC^n)$ spanned by the image of $\delta$ is uniquely decomposed as a finite sum $\alpha = \sum_i \theta_i \cdot \delta_{x_i}$.
We define a linear map $\iota : E \to \calG(\bC^n)$ by the equation $\iota(\alpha) = \sum_i \theta_i \cdot \lseg x_i \rseg - \left( \sum_i \theta_i \right) \cdot \lseg 0 \rseg = \partial \sum_i \theta_i \cdot \lseg 0,x_i \rseg$.
We shall refer to part (A) of \ref{3.2} to establish the metric claim in {\bf (ii)} and {\bf (iii)} below.
\par 
{\bf (ii)}
Let $\vphi \in \calD^0(\Rn)$ be so that $\bM(d\vphi) \leq 1$ and let $\psi$ be the restriction to $\bC^n$ of $x \mapsto \vphi(x) - \vphi(o)$.
Notice that $\psi \in \rmLip_o(\bC^n)$ and $\rmLip \psi \leq 1$.
Moreover, $\la \vphi , \iota(\alpha) \ra = \sum_i \theta_i \cdot \psi(x_i) = \alpha(\psi) \leq \|\alpha\|_L^*$.
Since $\vphi$ is arbitrary, we conclude that $\bG(\iota(\alpha)) \leq \|\alpha\|_L^*$.
\par 
{\bf (iii)}
Let $u \in \rmLip_o(\bC^n)$ be such that $\rmLip u \leq 1$.
Choose an extension $\hat{u}$ of $u$ to $\Rn$ so that $\rmLip \hat{u} \leq 1$ and $\hat{u}$ has compact support.
Abbreviate $\vphi_\veps = \phi_\veps * u$.
Thus, $\vphi_\veps \in \calD^0(\Rn)$ and $\bM(d\vphi_\veps) \leq 1$.
Note that $\la u , \alpha \ra = \sum_i \theta_i \cdot (u(x_i) - u(o)) = \lim_{\veps \to 0^+} \sum_i \theta_i \cdot ( \vphi_\veps(x_i) - \vphi_\veps(o)) = \lim_{\veps \to 0^+} \la \vphi_\veps , \iota(\alpha) \ra \leq \bG(\iota(\alpha))$.
Since $u$ is arbitrary, $\|\alpha\|_L^* \leq \bG(\iota(\alpha))$.
\par 
{\bf (iv)} 
It follows from {\bf (ii)} and {\bf (iii)} that $\iota$ is an isometry.
Since $E$ is dense in $\calF(\bC^n)$ (by definition of the latter) and $\calG(\bC^n)$ is Banach (by \ref{5.2}(A)), we infer that $\iota$ extends to a linear isometry $\hat{\iota} : \calF(\bC^n) \to \calG(\bC^n)$.
Therefore, it remains to prove that the image of $\hat{\iota}$ is $\bG$-dense.
As $\iota(E) \subset \rmim \hat{\iota}$, this is an immediate consequence of \ref{6.2}.
\end{proof}

\section{Qualitative $\bG$-compactness by means of a function $\bkappa$}

\begin{Empty}[Definition of $\bkappa(T,\veps)$]
\label{8.1}
Based on \ref{4.4} we make the following definition.
Let $\veps > 0$ and $T \in \calG(\bC^n)$ and define
\begin{multline*}
\bkappa(T,\veps) = \inf \big\{ \bM(\partial S) : S \in \bN_1(\Rn),\, \rmspt T \subset \bC^n, \,  \bG(T-\partial S) < \veps, \\ 
\text{ and } \bM(S) < \veps + \bG(T) \big\}.
\end{multline*}
This measures how close $T$ is of having finite mass at $\bG$-scale $\veps$, see \ref{8.0}(D).
\end{Empty}

\begin{Theorem}
\label{8.0}
The following hold.
\begin{enumerate}
\item[(A)] $(\forall T \in \calG(\bC^n))(\forall \veps > 0) : \bkappa(T,\veps) < \infty$.
\item[(B)] $(\forall T \in \calG(\bC^n)): \veps \mapsto \bkappa(T,\veps)$ is non-increasing.
\item[(C)] $(\forall T \in \calG(\bC^n)):$ the following are equivalent:
\begin{enumerate}
\item[(a)] $\bM(T) < \infty$;
\item[(b)] $(\forall \veps > 0): \bkappa(T,\veps) \leq \bM(T) < \infty$;
\item[(c)] $\lim_{\veps \to 0^+} \bkappa(T,\veps) < \infty$.
\end{enumerate}
\item[(D)] $(\forall T \in \calG(\bC^n)): \bkappa(T,\veps) = \inf \{ \bM(\hat{T}) : \hat{T} \in \calG(\bC^n) \text{ and } \bG(T-\hat{T}) < \veps \}$.
\item[(E)] $(\forall T,\hat{T} \in \calG(\bC^n))(\forall \veps, \hat{\veps} \in \R^+_0) : \bG(T-\hat{T}) < \veps \Rightarrow \bkappa(T,\veps+\hat{\veps}) \leq \bkappa(\hat{T},\hat{\veps})$.
\item[(F)] $(\forall T_1,T_2 \in \calG(\bC^n))(\forall \veps_1, \veps_2 \in \R^+_0) : \bkappa(T_1 + T_2, \veps_1 + \veps_2) \leq \bkappa(T_1,\veps_1) + \bkappa(T_2,\veps_2)$.
\item[(G)] $(\forall T \in \calG(\bC^n))(\forall \veps > 0)(\forall \lambda \in \R_0): \bkappa(\lambda \cdot T , |\lambda| \cdot \veps) = |\lambda| \cdot \bkappa(T,\veps)$.
\item[(H)] $(\forall T \in \calG(\bC^n)): \veps \mapsto \bkappa(T,\veps)$ is convex.
\item[(I)] $(\forall T \in \calG(\bC^n))(\forall \veps > 0): \bkappa(T,\veps) = 0 \iff \bG(T) \leq \veps$.
\end{enumerate}
\end{Theorem}

\begin{proof}
{\bf (i)}
Conclusion (A) is an immediate consequence of \ref{4.4} whereas conclusion (B) readily follows from the definition.
\par
{\bf (ii)}
Here, we show that $(a) \Rightarrow (b)$ in conclusion (C).
Assume that $T \in \calG(\bC^n)$ and $\bM(T) < \infty$.
Given $\veps > 0$ there exists $S \in \bM\bF_1(\Rn)$ such that $\partial S = T$, $\rmspt S \subset \bC^n$, and $\bM(S) < \veps + \bG^\times(T) = \bG(T)$, by \ref{5.2}(C).
Since $\bM(\partial S) = \bM(T) < \infty$, it follows that $S \in \bN_1(\Rn)$ and $S$ is, therefore, a competitor in the definition of $\bkappa(T,\veps)$.
Hence, $\bkappa(T,\veps) \leq \bM(\partial S) = \bM(T)$.
\par 
{\bf (iii)}
Since $(b) \Rightarrow (c)$ is trivial, in order to complete the proof of conclusion (C) it remains only to show that $(c) \Rightarrow (a)$.
Abbreviate $\gamma = \lim_{\veps \to 0^+} \bkappa(T,\veps)$.
For each $j=1,2,\cdots$ let $S_j \in \bN_1(\Rn)$ be a competitor in the definition of $\bkappa(T,1/j)$ such that $\bM(\partial S_j) \leq 1 + \gamma$.
Note that $\partial S_j \to T$ weakly* as $j \to \infty$, thus, $\bM(T) \leq \liminf_j \bM(\partial S_j) \leq 1 + \gamma$.
\par 
{\bf (iv)}
We turn to proving conclusion (D) and, to this end, we let $\hat{\bkappa}(T,\veps)$ be the quantity on the right side of the equality that we ought to establish.
If $S$ is a competitor in the definition of $\bkappa(T,\veps)$ then $\hat{T}=\partial S$ is clearly a competitor in the definition of $\hat{\bkappa}(T,\veps)$.
Therefore, $\bM(\partial S)=\bM(\hat{T}) \geq \hat{\bkappa}(T,\veps)$.
Since $S$ is arbitrary, we infer that $\bkappa(T,\veps) \geq \hat{\bkappa}(T,\veps)$.
\par 
{\bf (v)}
In order to establish the reverse inequality, we let $\hat{T}$ be a competitor in the definition of $\hat{\bkappa}(T,\veps)$.
We choose $\hat{\veps}_0 > 0$ such that $\bG(T-\hat{T}) + \hat{\veps}_0 < \veps$ and we let $0 < \hat{\veps} < \hat{\veps}_0$.
In turn, referring to the definition of $\bkappa(\hat{T},\hat{\veps})$ we choose $S \in \bN_1(\Rn)$ such that $\rmspt S \subset \bC^n$, $\bG(\hat{T} - \partial S) < \hat{\veps}$, $\bM(S) < \hat{\veps} + \bG(\hat{T})$, and $\bM(\partial S) < \hat{\veps} + \bkappa(\hat{T},\hat{\veps})$.
Observe that $\bG(T - \partial S) \leq \bG(T-\hat{T}) + \bG(\hat{T} - \partial S) < \bG(T - \hat{T}) + \hat{\veps} < \veps$ and $\bM(S) < \hat{\veps} + \bG(\hat{T}) \leq \hat{\veps} + \bG(\hat{T} - T) + \bG(T) < \veps + \bG(T)$, whence, $S$ is in competition for the definition of $\bkappa(T,\veps)$.
Accordingly, $\bkappa(T,\veps) \leq \bM(\partial S) < \hat{\veps} + \bkappa(\hat{T},\hat{\veps}) \leq \hat{\veps} + \bM(\hat{T})$, by (C)(b).
Since $\hat{\veps}$ is arbitrary, we infer that $\bkappa(T,\veps) \leq \bM(\hat{T})$.
Finally, $\bkappa(T,\veps) \leq \hat{\bkappa}(T,\veps)$, by the arbitrariness of $\hat{T}$.
\par 
{\bf (vi)}
In order to prove conclusion (E), we use the characterization of $\bkappa$ established in (D).
Let $\check{T} \in \calG(\bC^n)$ be such that $\bG(\hat{T} - \check{T}) < \hat{\veps}$.
Thus, $\bG(T - \check{T}) \leq \bG(T-\hat{T}) + \bG(\hat{T} - \check{T}) < \veps + \hat{\veps}$ and it follows that $\bkappa(T,\veps + \hat{\veps}) \leq \bM(\check{T})$.
Since $\check{T}$ is arbitrary, the proof is complete.
\par 
{\bf (vii)}
We refer again to the characterization of $\bkappa$ established in (D).
Let $\hat{T}_i \in \calG(\bC^n)$ be such that $\bG(T_i - \hat{T}) < \veps_i$, $i=1,2$.
Then $\bG[(T_1+T_2) - (\hat{T}_1 + \hat{T}_2)] \leq \bG(T_1 - \hat{T}_1) + \bG(T_2 - \hat{T}_2) < \veps_1 + \veps_2$.
Hence, $\bkappa(T_1+T_2,\veps_1+\veps_2) \leq \bM(\hat{T}_1 + \hat{T}_2) \leq \bM(\hat{T}_1) + \bM(\hat{T}_2)$.
By the arbitrariness of $\hat{T}_1$ and $\hat{T}_2$, we conclude that (F) holds.
\par 
{\bf (viii)}
Conclusion (G) rests on the simple observation that for all $T \in \calG(\bC^n)$, $\bG(T-\hat{T}) < \veps$ if and only if $\bG(\lambda \cdot T - \lambda \cdot \hat{T}) < |\lambda| \cdot \veps$ and $\bM(\lambda \cdot \hat{T}) = |\lambda| \cdot \bM(\hat{T})$.
\par 
{\bf (ix)}
We prove (H).
Let $T \in \calG(\bC^n)$, $0 < \veps_1 < \veps_2$, and $0 < \lambda < 1$.
It follows from (F) and (G) that
\begin{multline*}
\bkappa(T,\lambda \cdot \veps_1 + (1-\lambda) \cdot \veps_2) = \bkappa( \lambda \cdot T + (1-\lambda) \cdot T, \lambda \cdot \veps_1 + (1-\lambda) \cdot \veps_2) \\ 
\leq \bkappa(\lambda \cdot T , \lambda \cdot \veps_1) + \bkappa( (1-\lambda) \cdot T , (1-\lambda) \cdot \veps_2)
 = \lambda \cdot \bkappa(T,\veps_1) + (1-\lambda) \cdot \bkappa(T,\veps_2).
\end{multline*}
\par 
{\bf (x)}
Let $T \in \calG(\bC^n)$ and $\veps > 0$.
If $\bG(T) < \veps$ then $\hat{T}=0$ is a competitor in the characterization (D) of $\bkappa(T,\veps)$, therefore, $\bkappa(T,\veps) \leq \bM(\hat{T}) = 0$.
We also have $\bkappa(T,\bG(T))=0$, since $\veps \mapsto \bkappa(T,\veps)$ is continuous, by (H).
Assume now that $\bkappa(T,\veps)=0$.
Choose a sequence $\la \hat{T}_j \ra_j$ in $\calG(\bC^n)$ such that $\bG(T - \hat{T}_j) < \veps$ and $\bM(\hat{T}_j) \to 0$ as $j \to \infty$.
Since $T-\hat{T}_j \to T$ weakly* as $j \to \infty$, we infer that $\bG(T) \leq \liminf_j \bG(T-\hat{T}_j) \leq \veps$.
This completes the proof of (I).
\end{proof}

\begin{Theorem}
\label{8.3}
Let $\calT \subset \calG(\bC^n)$ be $\bG$-closed.
The following conditions are equivalent.
\begin{enumerate}
\item[(A)] $\calT$ is $\bG$-compact.
\item[(B)] There exist a positive real number $\Gamma$ and a non-increasing function $\kappa : \R^+_0 \to \R^+$ with the following properties.
\begin{enumerate}
\item[(a)] $(\forall T \in \calT) : \bG(T) \leq \Gamma$.
\item[(b)] $(\forall T \in \calT)(\forall \veps > 0) : \bkappa(T,\veps) \leq \kappa(\veps)$.
\end{enumerate}
\end{enumerate}
\end{Theorem}

\begin{proof}[Proof that $(A) \Rightarrow (B)$]
Condition (a) is obvious.
Let $\la\veps_j\ra_{j \geq 1}$ be a decreasing null sequence of real numbers and $\veps_0 = \infty$.
Owing to the compactness of $\calT$, for each $j \geq 1$ there exists a finite set $\calT_j \subset \calT$ with the following property: For all $T \in \calT$ there exists $\hat{T} \in \calT_j$ such that $\bG(T-\hat{T}) < \veps_j$.
In particular, $\bkappa(T,2\cdot \veps_j) \leq \bkappa(\hat{T},\veps_j)$, according to \ref{8.0}(E).
We define $\kappa_j = \max \{ \bkappa(\hat{T},\veps_j) : \hat{T} \in \calT_j \}$.
As there is no restriction to assume that $\la \calT_j \ra_{j \geq 1}$ is non-decreasing, so is $\la \kappa_j \ra_{j \geq 1}$.
We now define $\kappa : \R^+_0 \to \R^+$ by requiring that $\kappa(r) = \kappa_j$ whenever $2 \cdot \veps_j \leq r < 2 \cdot \veps_{j-1}$, $j \geq 1$.
In view of \ref{8.0}(B), it is now clear that $\kappa$ satisfies condition (b).
\end{proof}

\begin{proof}[Proof that $(B) \Rightarrow (A)$]
Let $\la T_j \ra_{j \geq 1}$ be a sequence in $\calT$ and let $\la \veps_k \ra_{k \geq 1}$ be a decreasing null sequence of real numbers.
For each $j$ and $k$ there exists $S_{j,k} \in \bN_1(\Rn)$ satisfying the following properties:
\begin{enumerate}
\item[(1)] $\rmspt S_{j,k} \subset \bC^n$;
\item[(2)] $\bG(T_j - \partial S_{j,k}) < \veps_k$;
\item[(3)] $\bM(S_{j,k}) < \veps_k + \bG(T_j) \leq \veps_1 + \Gamma$;
\item[(4)] $\bM(\partial S_{j,k}) < 1 + \bkappa(T_j , \veps_k) \leq 1 + \kappa(\veps_k)$.
\end{enumerate}
Thus, for each $k$ the sequence $ \la S_{j,k} \ra_j$ has uniformly bounded normal mass, by (3) and (4), and consists currents all supported in the same compact set, by (1).
Accordingly, the Federer-Fleming compactness theorem \cite[4.2.17(1)]{GMT} yields by induction on $k$ increasing sequences of integers $\la j_k(l) \ra_l$ and $S_k \in \bN_1(\Rn)$ such that $\la j_k(l) \ra_l$ is a subsequence of $\la j_{k-1}(l) \ra$ and $\bF(S_{j_k(l)} - S_k) \to 0$ as $l \to \infty$.
For each $k$ we choose an integer $l_0(k)$ such that $\bF(S_{j_k(l)} - S_k) < \veps_k$ whenever $l \geq l_0(k)$.
Finally, we define by induction on $k$ an increasing sequence of integers $\la j(k)\ra_k$ so that $j(k) = j_k(l_k)$ for some $l_k \geq l_0(k)$ and we check that it is a subsequence of $\la j_k(l)\ra_l$ for each $k$.
We now observe that for all $k_1 < k_2$ we have, in view of (2) and \ref{3.1}(B),
\begin{multline*}
\bG \left( T_{j(k_1)} - T_{j(k_2)} \right) \\ \leq
\bG \left( T_{j(k_1)} - \partial S_{j(k_1),k_1} \right) + \bG \left( \partial S_{j(k_1),k_1} - \partial S_{j(k_2),k_1} \right) + \bG \left( \partial S_{j(k_2),k_1} -  T_{j(k_2)} \right) \\
< \veps_{k_1} + \bF \left(  S_{j(k_1),k_1} -  S_{j(k_2),k_1} \right) + \veps_{k_1} < 4 \cdot \veps_{k_1}.
\end{multline*}
Thus, $\la T_{j(k)} \ra_k$ is $\bG$-Cauchy, whence, $\bG$-convergent, by \ref{5.2}(A).
As $\la T_j \ra_j$ is arbitrary, the proof is complete.
\end{proof}

\begin{Empty}[$(\Gamma,\kappa)$-boundedness]
\label{8.5}
Let $\Gamma > 0$ and $\kappa : \R^+_0 \to \R^+$ be a non-increasing function.
We say that $\calT \subset \calG(\bC^n)$ is {\em $(\Gamma,\kappa)$-bounded} provided it satisfies conditions (a) and (b) of \ref{8.3}(B).
In other words, the $(\Gamma,\kappa)$-bounded subsets of $\calG(\bC^n)$, corresponding to all choices of $\Gamma$ and $\kappa$, are precisely those which are relatively $\bG$-compact.
\end{Empty}

\begin{Empty}[Modulus of smoothness]
For $h \in \Rn$ we abbreviate $\btau_h : \Rn \to \Rn : x \mapsto x+h$.
Let $\eta : \Rn \to \bigwedge_1\Rn$ be and Lebesgue-summable.
For each $r > 0$ define
\begin{equation*}
\rmosc_1(\eta,r) = \sup \left\{ \int_{\Rn} |\eta - \eta \circ \btau_h | d\calL^n : h \in \Rn \text{ and } |h| \leq r \right\}.
\end{equation*}
The function $\R^+_0 \to \R^+ : r \mapsto \rmosc_1(\eta,r)$ is called the {\em modulus of smoothness in $L_1$ of $\eta$}.
\par 
If $\eta$ is as in the previous paragraph and $\omega(r) = \rmosc_1(\eta,r)$ then $\omega$ satisfies the following properties.
It is non-decreasing, subadditive (\ie $\omega(r+s) \leq \omega(r) + \omega(s)$ for all $0 < r < r+s$), and $\lim_{r \to 0^+} \omega(r) = 0$.
With regard to the next result it is useful to note that a modulus of smoothness in $L_1$ does not need to satisfy any other property in a neighborhood of $0$, in particular, its decay near $0$ can be arbitrarily slow.
Specifically, if $\omega$ satisfies the three properties stated above then there exists a Lebesgue-summable $\eta : \R \to \R$ such that $\omega(r) \leq \rmosc_1(\eta,r)$ for all $0 < r \leq 1$, see \cite[Theorem 3 p.406]{KOL.75}.
\end{Empty}

\begin{Theorem}
\label{8.4}
There exists $\bc_{\theTheorem}(n) > 0$ with the following property.
Let $T \in \calG(\bC^n)$ and let $\eta$ be associated with $T$ and $\veps=1$ in \ref{4.2}.
Then for every $\veps > 0$ we have
\begin{equation*}
\bkappa(T,\rmosc_1(\eta,\veps)) \leq \frac{\bc_{\theTheorem}(n)\cdot(1+ \bG(T))}{\veps}. 
\end{equation*}
\end{Theorem}

\begin{proof}
{\bf (i)}
Abbreviate $S = \calL^n \wedge \eta$.
Define $S_\veps = \phi_\veps * S$.
Here, we aim to find an upper bound for $\bM(S - S_\veps)$.
Notice that $\la \omega , S_\veps \ra = \la \phi_\veps * \omega , S \ra = \int_{\Rn} \la \phi_\veps * \omega , \eta \ra d\calL^n = \int_{\Rn} \la \omega , \phi_\veps * \eta \ra d\calL^n$ for all $\omega \in \calD^1(\Rn)$, therefore, $S_\veps = \calL^n \wedge (\phi_\veps * \eta)$.
It follows that
\begin{multline*}
\bM(S - S_\veps) = \int_{\Rn} |\eta - \phi_\veps * \eta|d\calL^n 
= \int_{\Rn} \left| \int_{\Rn} \phi_\veps(h)\cdot \big( \eta(x) - \eta(x-h) \big) \right|d\calL^n(x) \\
\leq \int_{\Rn} \left( \int_{\Rn} | \eta(x)-\eta(x-h)| d\calL^n(x)\right) \phi_\veps(h) d\calL^n(h) \leq \rmosc_1(\eta,\veps).
\end{multline*}
\par 
{\bf (ii)}
Notice that $S_\veps$ has compact support, because so does $S$, and is normal.
Thus, $\pi_\# S \in \bN_1(\Rn)$ and $\rmspt \pi_\# S \subset \bC^n$.
We define $\hat{T} = \partial \pi_\# S_\veps$ and we observe that $\hat{T} \in \calG(\bC^n)$.
Furthermore, since $T = \pi_\# T$, we have
\begin{equation*}
\bG(T - \hat{T}) = \bG(\partial \pi_\#(S-S_\veps)) \leq \bM(\pi_\#(S-S_\veps)) \leq \bM(S-S_\veps) \leq \rmosc_1(\eta,\veps),
\end{equation*}
by {\bf (i)}.
\par 
{\bf (iii)}
Next, we determine an upper bound for $\bM(\partial S_\veps)$.
If $\vphi \in \calD^0(\Rn)$ then $\la \vphi , \partial S_\veps \ra = \la \phi_\veps * d\vphi , S \ra = \la (d\phi_\veps)*\vphi , S \ra \leq \bM[(d\phi_\veps)*\vphi] \cdot \bM(S)$.
Furthermore, for each $x \in \Rn$ we have
\begin{multline*}
[(d\phi_\veps)*\vphi](x) = \int_{\Rn} (d\phi_\veps)(h)\cdot \vphi(x-h) d\calL^n(h) \leq \bM(\vphi) \int_{\Rn} |d\phi_\veps(h)|d\calL^n(h)
\\ = \bM(\vphi) \int_{\Rn} |d(\phi \circ \bmu_{\veps^{-1}})(h)|\veps^{-n}d\calL^n(h)  = \frac{\bM(\vphi)}{\veps} \int_{\Rn} |[(d\phi) \circ \bmu_{\veps^{-1}}](h)| \veps^{-n} d\calL^n(h) \\
= \frac{\bM(\vphi)}{\veps} \int_{\Rn} |d\phi|d\calL^n.
\end{multline*}
It follows that 
\begin{equation*}
\bM(\partial S_\veps) \leq \frac{\bc_{\ref{8.4}}(n) \cdot \bM(S)}{\veps}
\end{equation*}
where $\bc_{\ref{8.4}}(n) = \int_{\Rn} |d\phi|d\calL^n$.
Recalling that $\bM(S) = \bM(\calL^n \wedge \eta) \leq 1 + \bG(T)$ we conclude that
\begin{equation*}
\bM(\hat{T}) = \bM(\pi_\# \partial S_\veps) \leq  \bM(\partial S_\veps) \leq \frac{\bc_{\ref{8.4}}(n) \cdot (1+ \bG(T))}{\veps}.
\end{equation*}
\par 
{\bf (iv)}
We infer from \ref{8.0}(D), {\bf (ii)} and {\bf (iii)} that $\bkappa(T,\lambda + \rmosc_1(\eta,\veps)) \leq \bM(\hat{T}) \leq \frac{\bc_{\ref{8.4}}(n) \cdot (1+ \bG(T))}{\veps}$ for all $\lambda > 0$.
As $\bkappa(T,\cdot)$ is continuous, according to \ref{8.0}(H), the conclusion ensues.
\end{proof}

\section{Quantitative $\bG$-compactness by means of the deformation theorem}

\begin{Empty}[Notations]
\label{9.1}
Throughout this section, we enforce the following notations.
\begin{enumerate}
\item[(A)] $k$ is a positive integer and $\veps_k = \frac{1}{k}$.
\item[(B)] $\Sigma^{n,0}_k = \bC^n \cap (\veps_k \cdot \Z^n)$. Observe that $\rmcard \Sigma^{n,0}_k = (2 k + 1)^n$.
\item[(C)] We let $\Sigma^{n,1}_k$ be the 1-skeleton of the decomposition of $[-1,1]^n$ into congruent cubes of edge length equal to $\frac{1}{k}$. In other words, $\Sigma^{n,1}_k$ consists of all line segments $\Rn \cap \{ a + t\cdot (b-a) : 0 \leq t \leq 1 \}$ corresponding to $a,b \in \Sigma^{n,0}_k$ that are {\it adjacent}, \ie $a$ and $b$ differ by exactly one coordinate and $|a-b| = \frac{1}{k}$. Notice that in this case $\rmspt [\partial ( \theta \cdot \lseg a , b \rseg )] \subset \Sigma^{n,0}_k$ whenever $\theta \in \R$. 
\end{enumerate}
\end{Empty}

\begin{Theorem}
\label{9.2}
There exists $\bc_{\theTheorem}(n) \geq 1$ such that for all $k \in \N_0$ and all member $P$ of the vector space
\begin{equation*}
\bP_0(\Rn) \cap \left\{ P : \rmspt P \subset \Sigma^{n,0}_k \text{ and } \chi(P)=0 \right\}
\end{equation*}
we have
\begin{equation*}
\frac{\bG(P)}{\sqrt{n}} \leq \bM(P) \leq \bc_{\theTheorem}(n) \cdot k^{2  n} \cdot \bG(P).
\end{equation*}
\end{Theorem}

\begin{proof}
{\bf (i)}
The inequality $\bG(P) \leq \sqrt{n} \cdot \bM(P)$ follows from \ref{3.3} applied with $a = 0$ and $S = \lseg 0 \rseg \cone P$.
Note that $\partial S = P$, since $\chi(P)=0$, and $\rmspt S \subset \bC^n$.
\par 
{\bf (ii)}
{\it Assume that $F \subset \R$ is a finite set, $\rho > 0$ is so that $|x-y| \geq \rho$ whenever $x,y \in F$ are distinct, and $P \in \bP_0(\R)$ satisfies $\rmspt P \subset F$ and $\chi(P)=0$.
Then
\begin{equation*}
\bM(P) \leq \left( \frac{\rmcard F}{\rho} \right) \cdot \bG(P).
\end{equation*}
}
\par 
This number is devoted to the proof of that claim.
Notice that if $F$ is empty or a singleton then $P=0$ and the conclusion is obvious.
We henceforth assume that $\rmcard F \geq 2$.
Choose a numbering of $F = \{x_1,\cdots,x_N\}$ so that $x_1 < \cdots < x_N$ and real numbers $\theta_1,\cdots,\theta_N$ so that $P = \sum_{i=1}^N \theta_i \cdot \lseg x_i \rseg$.
\par 
We will now show that there exists a unique compactly supported $S \in \calD_1(\R)$ such that $P = \partial S$.
Uniqueness follows from the constancy theorem \cite[4.1.7 p.357]{GMT}: If $S_1$ and $S_2$ are compactly supported and have boundary $P$ then $\partial(S_1-S_2)=0$ and $S_1-S_2$ is compactly supported, hence, $S_1-S_2=0$.
We now show that $S = \sum_{i=1}^{N-1} \gamma_i \cdot \lseg x_i , x_{i+1} \rseg$ for some real numbers $\gamma_1,\cdots,\gamma_{N-1}$.
As $\partial S = \sum_{i=1}^{N-1} ( \gamma_i \cdot \lseg x_{i+1} \rseg - \gamma_i \cdot \lseg x_i \rseg ) = (-\gamma_1) \cdot \lseg x_1 \rseg + \sum_{i=2}^{N-1} (\gamma_{i-1} - \gamma_i) \cdot \lseg x_i \rseg + \gamma_{N-1} \cdot \lseg x_N \rseg$, we see that the requirement $\partial S = P$ is equivalent to the equations $-\gamma_1 = \theta_1$, $\gamma_{i-1} - \gamma_i = \theta_i$ for all $i \in \{2,\cdots,N-1\}$, and $\gamma_{N-1} = \theta_N$.
The first $N-1$ equations determine the values of
\begin{equation}
\label{eq.4}
\gamma_i = - \sum_{j=1}^i \theta_j
\end{equation}
and the last one is satisfied, since $\chi(P)=0$.
\par 
We infer from the previous paragraph that $\bG(P) = \bM(S)$.
Moreover,
\begin{equation}
\label{eq.5}
\bM(S) = \sum_{i=1}^{N-1} |\gamma_i| \cdot |x_{i+1}-x_i | \geq \rho \cdot \sum_{i=1}^{N-1} \left| \sum_{j=1}^i \theta_j \right|,
\end{equation}
in view of \eqref{eq.4}.
\par 
Let $i_0 \in \{1,\cdots,N\}$ be so that $|\theta_{i_0}| = \max \{ | \theta_i | : i=1,\cdots,N \}$.
We claim that $ \sum_{i=1}^{N-1} \left| \sum_{j=1}^i \theta_j \right| \geq | \theta_{i_0}|$.
If $i_0=1$ then this is obvious. 
If not, then $\left| \sum_{j=1}^{i_0} \theta_j \right| \geq |\theta_{i_0}| - \left| \sum_{j=1}^{i_0-1} \theta_j \right|$ and our claim ensues.
\par 
It follows from \eqref{eq.5} that $\bM(S) \geq \rho \cdot | \theta_{i_0} |$.
Since $\bM(P) = \sum_{i=1}^N |\theta_i| \leq N \cdot |\theta_{i_0}|$, the proof of the claim {\bf (ii)} is complete.
\par 
Applying {\bf (ii)} to the case when $F = \Sigma^{1,0}_k$ and $\rho = \veps_k$ we infer that $\bM(P) \leq k  (2 k + 1) \cdot \bG(P)$.
Whence the second inequality in the statement of the theorem holds for $n=1$ with $\bc_{\ref{9.2}}(1) = 3$.
\par 
{\bf (iii)}
{\it Assuming that $n \geq 2$, there exists a unit vector $u_k \in \Rn$ such that for all distinct $x , y \in \Sigma^{n,0}_k$ we have 
\begin{equation*}
|(x-y) \ip u_k | \geq  \frac{\veps_k^{n}}{2 \cdot \sqrt{2} \cdot  5^{n-1}}.
\end{equation*} 
}
\par 
This number is devoted to the proof of that claim.
We define $\alpha_k = \frac{\veps_k}{4 + \veps_k}$ and we note that $0 < \alpha_k \leq \frac{1}{2}$ and $\frac{\alpha_k}{1-\alpha_k} = \frac{\veps_k}{4}$.
Next, we define $s_j = \alpha_k^{j-1}$ for $j \in \{1,\cdots,n\}$, $v_k = \sum_{j=1}^n s_j \cdot e_j$ and $u_k = \frac{v_k}{|v_k|}$.
Notice that 
\begin{equation}
\label{eq.3}
s_n \geq \left( \frac{\veps_k}{5} \right)^{n-1} \text{ and }|v_k| \leq \sqrt{2}.
\end{equation}
Let $x,y \in \Sigma^{n,0}_k$ be such that $x \neq y$.
Decompose $x - y = \sum_{j=1}^n t_j \cdot e_j$.
For all $j$ observe that $|t_j| \leq 2$ and either $t_j = 0$ or $|t_j| \geq \veps_k$.
Since $x \neq y$, there are $j \in \{1,\cdots,n\}$ such that $t_j \neq 0$. 
Let $j_0$ be the smallest of them.
Thus, $(x-y) \ip v_k = \sum_{j=j_0}^n t_j \cdot s_j$.
If $j_0=n$ then $|(x-y) \ip v_k| = |t_n \cdot s_n| \geq \veps_k \cdot s_n$ and the claim holds, according to \eqref{eq.3}.
We henceforth assume that $j_0 < n$.
In that case,
\begin{multline*}
| (x-y) \ip v_k | = \left| \sum_{j=j_0}^n t_j \cdot s_j \right| \geq |t_{j_0}| \cdot s_{j_0} -  \sum_{j=j_0+1}^n |t_j| \cdot s_j \geq \veps_k \cdot \alpha_k^{j_0-1} - 2 \cdot \sum_{j=j_0+1}^n \alpha_k^{j-1} \\
\geq \veps_k \cdot \alpha_k^{j_0-1} - 2 \cdot \alpha_k^{j_0-1} \cdot \left( \frac{\alpha_k}{1-\alpha_k}\right) = \frac{\veps_k}{2} \cdot s_{j_0}  \geq \frac{\veps_k}{2} \cdot s_n \geq   \frac{\veps_k^{n}}{2 \cdot 5^{n-1}} \\
\geq |v_k| \cdot \left( \frac{\veps_k^{n}}{2 \cdot \sqrt{2} \cdot  5^{n-1} }\right),
\end{multline*} 
in view of \eqref{eq.3}.
This completes the proof of {\bf (iii)}.
\par 
{\bf (iv)}
Assume $n \geq 2$ and choose $u_k$ according to {\bf (iii)}.
Define a smooth map $f : \Rn \to \R : x \mapsto \la x , u_k \ra$ and notice that $\rmLip f \leq 1$, since $u_k$ is a unit vector.
We now consider the distribution $f_\#P \in \bP_0(\Rn)$ in $\R$ and we notice that $\chi(f_\#P)=0$.
\par 
Observe that $f|_{\Sigma^{n,0}_k}$ is injective, by {\bf (iii)}.
As a first consequence, note that $\bM(f_\#P) = \bM(P)$, since $\rmspt P \subset \Sigma^{n,0}_k$.
A second consequence is that $\rmcard F = (2k+1)^n$, where we have abbreviated $F = f(\Sigma^{n,0}_k)$.
\par 
Letting $\rho = \frac{\veps_k^{n}}{2 \cdot \sqrt{2} \cdot  5^{n-1}}$, we note that if $x,y \in F$ are distinct then $x = f(\tilde{x})$ and $y = f(\tilde{y})$ for some distinct $\tilde{x},\tilde{y} \in \Sigma^{n,0}_k$ and, therefore, $|x-y| = |f(\tilde{x}) - f(\tilde{y})| = | \la \tilde{x} - \tilde{y} , u_k \ra | \geq \rho$, by {\bf (iii)}.
\par 
Finally,
\begin{equation*}
\bM(P) = \bM(f_\#P) \leq \left( \frac{\rmcard F}{\rho} \right) \cdot \bG(f_\#P) \leq 2 \cdot \sqrt{2} \cdot 5^{n-1} \cdot k^n(2k+1)^n \cdot \bG(P),
\end{equation*}
according to {\bf (iii)} and $\bG(f_\#P) \leq \bG(P)$.
Thus, the theorem holds in case $n \geq 2$ with \eg $\bc_{\ref{9.2}}(n) = 2 \sqrt{2} \cdot 5^{n-1} \cdot 3^n$.
\end{proof}

\begin{Empty}[The class $\calP^n_{k,\veps}$]
\label{9.3}
The class of distributions that we are about to define depends on $n \in \N_0$, the dimension of the ambient space $\Rn$, on $k \in \N_0$, as it pertains to the edge length of a subdivision of the unit cube of $\Rn$ (recall \ref{9.1} and \ref{9.2}), and on $0 < \veps \leq 1$.
\par 
We define a positive real number
\begin{equation}
\label{eq.10}
\hat{\veps}(n,k,\veps) = \frac{\veps}{2n \cdot (2k+1)^n}.
\end{equation}
Next, we define a subset of $\R$ as follows
\begin{equation}
\label{eq.11}
Z(n,k,\veps) =  \hat{\veps}(n,k,\veps) \cdot \Z .
\end{equation}
\par 
We are now ready to define $\calP^n_{k,\veps}$ by means of the following formula.
\begin{multline}
\label{eq.12}
\calP^n_{k,\veps} = \bP_0(\Rn) \cap \bigg\{ P :  P = \sum_{x \in \Sigma^{n,0}_k} \theta(x) \cdot \lseg x \rseg \text{ where } \theta(x) \in Z(n,k,\veps) \text{ for all } x \in \Sigma^{n,0}_k , \\ \chi(P) = 0, \text{ and } \bM(P) \leq k \cdot \veps \bigg\}
\end{multline}
\end{Empty}

\begin{Theorem}
\label{9.4}
There exists $0 < \bc_{\theTheorem}(n) < 1$ such that for all $k \in \N_0$ and $0 < \veps \leq 1$ the following hold.
\begin{enumerate}
\item[(A)] For all $P \in \calP^n_{k,\veps}$ we have
\begin{equation*}
\bM(P) \leq k \cdot \veps.
\end{equation*}
\item[(B)] For all $P \in \calP^n_{k,\veps}$ we have
\begin{equation*}
\bG(P) \leq \sqrt{n} \cdot k \cdot \veps.
\end{equation*}
\item[(C)] For all $P_1,P_2 \in \calP^n_{k,\veps}$ if $P_1 \neq P_2$ then
\begin{equation*}
\bG(P_1 - P_2) \geq \frac{\veps \cdot \bc_{\theTheorem}(n)}{k^{3n}}.
\end{equation*}
\end{enumerate}
\end{Theorem}

\begin{proof}
{\bf (i)}
Conclusion (A) is a direct consequence of the definition of $\calP^n_{k,\veps}$ whereas conclusion (B) follows from (A) and the inequality on the left side of \ref{9.2}.
\par 
{\bf (ii)}
In order to establish conclusion (C) we consider $P_i \in \calP^n_{k,\veps}$, $i=1,2$.
Thus, $P_i = \sum_{x \in \Sigma^n_{k,0}} \theta_i(x) \cdot \lseg x \rseg$ for some $\theta_i : \Sigma^{n,0}_{k} \to Z(n,k,\veps)$.
Notice that $\bM(P_1 - P_2) = \sum_{x \in \Sigma^{n,0}_{k}} |\theta_1(x) - \theta_2(x)|$.
If $P_1 \neq P_2$ then there exists $x_0 \in \Sigma^{n,0}_{k}$ such that $\theta_1(x_0) \neq \theta_1(x_0)$, whence, $|\theta_1(x_0) - \theta_2(x_0)| \geq \hat{\veps}(n,k,\veps)$, by \eqref{eq.11}.
Accordingly, $\bM(P_1 - P_2) \geq \hat{\veps}(n,k,\veps)$.
We then infer from the inequality on the right side of \ref{9.2} that $\bG(P_1-P_2) \geq \frac{\hat{\veps}(n,k,\veps)}{\bc_{\ref{9.2}}(n) \cdot k^{2n}}$.
It is now plain that $\bc_{\ref{9.4}}(n) = [2n \cdot 3^n \cdot \bc_{\ref{9.2}}(n)]^{-1}$ suits our needs.
\end{proof}

\begin{Theorem}
\label{9.5}
There exists $\bc_{\theTheorem}(n) \geq 1$ such that for all integer $k \geq \bc_{\theTheorem}(n)$, all $0 < \veps \leq 1$, and all $T \in \calG(\bC^n)$, among the following three conditions, the implications $(A) \Rightarrow (B) \Rightarrow (C)$ hold.
\begin{enumerate}
\item[(A)] We have
$
\bc_{\theTheorem}(n) \cdot \left( \bG(T) + \bkappa \left( T,  \veps\right) \right) < k \cdot\veps.
$
\item[(B)] There exists $P \in \calP^n_{k,\veps}$ such that $\bG(T-P) < 3 \cdot \veps$.
\item[(C)] We have
$
\left(\sqrt{n} + 4\right) \cdot \left( \bG(T) + \bkappa \left( T,  3 \cdot \veps\right) \right) < k \cdot\veps.
$
\end{enumerate}
\end{Theorem}

\begin{proof}[Proof that $(A) \Rightarrow (B)$]
{\bf (i)} {\it Choice of $\bc(n)$, $\bc_{\ref{9.5}}(n)$, and $\beta > 0$.}
We abbreviate $\bc(n) = 2n^4$ (which is the constant, called $\gamma$ there, appearing in the deformation theorem \cite[4.2.9]{GMT}, corresponding to the case $m=1$) and $\bc_{\ref{9.5}}(n)=2 \cdot \bc(n)$.
In the remaining part of this proof we assume that $T$, $k$ and $\veps$ satisfy condition (A) and that $k \geq \bc_{\ref{9.5}}(n)$.
Furthermore, we let $\beta > 0$ be so small that
\begin{equation*}
\bc(n) \cdot \left( \beta + \bG(T) + \bkappa \left( T,  \veps\right) \right) < k \cdot \left( \frac{\veps}{2} \right).
\end{equation*}
Since $k \geq 2 \cdot \bc(n)$, we have $\bc(n) \cdot \veps \leq k \cdot \left( \frac{\veps}{2} \right)$.
Thus,
\begin{equation*}
\bc(n) \cdot \left( \beta + \veps + \bG(T) + \bkappa \left( T,  \veps\right) \right) < k \cdot \veps.
\end{equation*}
\par 
{\bf (ii)}
{\it Choice of $S$ given $T$.}
There exists $S \in \bN_1(\Rn)$ such that $\rmspt S \subset \bC^n$ and the following hold.
\begin{enumerate}
\item[(1)] $\bG(T - \partial S) < \veps$;
\item[(2)] $\bM(S) < \veps + \bG(T)$;
\item[(3)] $\bM(\partial S) < \beta + \bkappa(T , \veps)$;
\item[(4)] $\bN(S) < \beta + \veps + \bG(T) + \bkappa(T,\veps)$. 
\end{enumerate}
The existence of such $S$ satisfying (1), (2), and (3) is a consequence of the definition of $\bkappa(T,\veps)$.
Inequality (4) is an immediate consequence of (2) and (3).
\par 
{\bf (iii)}
{\it Deforming $S$.}
Applying the deformation theorem \cite[4.2.9]{GMT} to $S \in \bN_1(\Rn)$ we obtain a $1$-dimensional (polyhedral) current $Q$ in $\Rn$ satisfying the following properties.
\begin{enumerate}
\item[(5)] $\rmspt Q \subset \bC^n$;
\item[(6)] $Q = \sum_{\sigma \in \Sigma^{n,1}_k} \theta_Q(\sigma) \cdot \lseg \sigma \rseg$ for some $\theta_Q : \Sigma^{n,1}_k \to \R$;
\item[(7)] $\bF(S-Q) \leq \veps_k \cdot \bc(n) \cdot \bN(S)$;
\item[(8)] $\bM(\partial Q) \leq \bc(n) \cdot \bM(\partial S)$.
\end{enumerate}
Specifically, we apply the deformation theorem with $m=1$, with $\veps_k$ in place of $\veps$, and with $S$ in place of $T$; our $Q$ here is the $P$ there.
\par 
With regard to (5), we note that the construction in the proof of \cite[4.2.9]{GMT} moves the centers of projections in each sub-cube rather than translating $S$ itself (as is done \eg in \cite{WHI.99.deformation}). 
\par 
Condition (6) is a rewriting of \cite[4.2.9(5)]{GMT}, recalling the form of the constancy theorem applied in \cite[4.2.5, middle of p.402]{GMT}.
Notice that our notation $\lseg \sigma \rseg$ is ambiguous, since it could mean two opposite currents according to the orientation we give to $\sigma$ (we haven't given any) but this ambiguity is absorbed in the sign of the coefficient $\theta_Q(\sigma)$.
At any rate, the only way we will use (6) is to infer that $\rmspt \partial S \subset \Sigma^{n,0}_k$ (recall \ref{9.1}(C)), which is already stated in \cite[4.2.9(5)]{GMT}.
\par 
Inequality (7) follows from the definition of the flat norm $\bF$, \cite[4.2.9(1)]{GMT} and the two inequalities in the bottom line of \cite[4.2.9(2)]{GMT} (this is where we need $\bc(n)$ to be equal to $\gamma$ there).
Similarly, inequality (8) follows from the inequality on the top of the right column of \cite[4.2.9(2)]{GMT}.
\par 
{\bf (iv)}
{\it Upper bound for $\bG(T - \partial Q)$.}
Here, we show that $\bG(T - \partial Q) < 2 \cdot \veps$.
Indeed,
\begin{equation*}
\begin{aligned}
\bG(T - \partial Q) & \leq \bG(T - \partial S) + \bG(\partial S - \partial Q) && \\
& < \veps + \bF(S - Q) && \text{(by {\bf (ii)}(1) and \ref{3.1}(B))} \\
&  \leq \veps + \veps_k \cdot \bc(n) \cdot \bN(S) && \text{(by {\bf (iii)}(7))} \\
& \leq \veps + \veps_k \cdot \bc(n) \cdot ( \beta + \veps + \bG(T) + \bkappa(T,\veps)) && \text{(by {\bf (ii)}(4))} \\
& \leq \veps + \veps &&\text{(by {\bf (i)})}.
\end{aligned}
\end{equation*}
\par 
{\bf (v)}
{\it Defining $P$ by modifying the multiplicities of $\partial Q$.}
We start by estimating from above the mass of $\partial Q$ as follows.
\begin{equation}
\label{eq.6}
\begin{aligned}
\bM(\partial Q) & \leq \bc(n) \cdot \bM(\partial S) && \text{(by {\bf (iii)}(8))} \\
& < \bc(n) \cdot \left( \beta + \bkappa (T,\veps)  \right) &&\text{(by {\bf (ii)}(3))} \\
& < k \cdot \left( \frac{\veps}{2}\right) &&\text{(by {\bf (i)})}.
\end{aligned}
\end{equation}
\par
Recall that 
\begin{equation*}
\partial Q = \sum_{x \in \Sigma^{n,0}_k} \theta_{\partial Q}(x) \cdot \lseg x \rseg
\end{equation*}
for some $\theta_{\partial Q} : \Sigma^{n,0}_k \to \R$ (this follows from {\bf (iii)}(6)).
For each $x \in \Sigma^{n,0}_{k}$ we choose $\theta_R(x) \in Z(n,k,\veps)$ such that $|\theta_{\partial Q}(x) - \theta_R(x)| < \hat{\veps}(n,k,\veps)$.
Now, we define $R \in \bP_0(\Rn)$ by
\begin{equation*}
R = \sum_{x \in \Sigma^{n,0}_k} \theta_{R}(x) \cdot \lseg x \rseg ,
\end{equation*}
Note that 
\begin{equation}
\label{eq.9}
\bM(R - \partial Q) = \sum_{x \in \Sigma^{n,0}_k} \left|\theta_{R}(x) - \theta_{\partial Q}(x) \right| \leq \hat{\veps}(n,k,\veps) \cdot \rmcard \Sigma^{n,0}_k = \frac{\veps}{2n},
\end{equation}
by \eqref{eq.10}.
Since $\chi(\partial Q ) = 0$, by \ref{P.AM}(C), we observe that
\begin{equation*}
|\chi(R)| = |\chi(R - \partial Q)| \leq \bM(R - \partial Q) \leq \frac{\veps}{2n},
\end{equation*}
in view of \ref{P.AM}(B).
Finally, we are ready to define $P$ by means of the formula
\begin{equation*}
P = R - \chi(R) \cdot \lseg x_0 \rseg
\end{equation*}
where $x_0 \in \Sigma^{n,0}_k$ is chosen arbitrarily.
\par 
{\bf (vi)} 
{\it Showing that $P \in \calP^n_{k,\veps}$.}
Note that $\rmspt P \subset \Sigma^{n,0}_k$, thus, $P = \sum_{x \in \Sigma^{n,0}_k} \theta_P(x) \cdot \lseg x \rseg$.
In fact, $\theta_P(x) = \theta_R(x)$ for all $x \neq x_0$ and $\theta_P(x_0) = \theta_R(x_0) - \chi(R)$.
Evidently, $\chi(R) \in Z(n,k,\veps)$, therefore, $\theta_P(x) \in Z(n,k,\veps)$ for all $x \in \Sigma^{n,0}_k$.
Furthermore, $\chi(P) = 0$, by definition of $P$.
Finally,
\begin{equation}
\label{eq.13}
\bM(P - \partial Q) \leq \bM(P - R) + \bM(R - \partial Q) \leq |\chi(R)| + \bM(R - \partial Q) \leq \frac{\veps}{n}
\end{equation}
and, in view of \eqref{eq.6},
\begin{equation*}
\bM(P) \leq \bM(\partial Q) + \frac{\veps}{n} < k \cdot \left( \frac{\veps}{2} \right) + \frac{\veps}{n} \leq k \cdot \veps,
\end{equation*}
where we used the inequality $\frac{1}{n} \leq \frac{k}{2}$ (recall that $k \geq \bc_{\ref{9.5}}(n) = 2 \cdot \bc(n) > 2$).
This completes the proof that $P \in \calP^n_{k,\veps}$.
\par 
{\bf (vii)}
{\it Upper bound for $\bG(T - P)$.}
Applying \ref{9.2} to $P - \partial Q$ we infer that $\bG(P - \partial Q) \leq \sqrt{n} \cdot \bM(P - \partial Q) \leq \veps$, in view of \eqref{eq.13}.
Finally, $\bG(T - P) \leq \bG(T - \partial Q) + \bG(\partial Q - P) < 3 \cdot \veps$, according to {\bf (iv)}.
\end{proof}

\begin{proof}[Proof that $(B) \Rightarrow (C)$]
Note that $\bG(T) \leq \bG(T-P) + \bG(P) < 3 \cdot \veps + \sqrt{n} \cdot k \cdot \veps \leq \left( \sqrt{n} + 3 \right) \cdot k \cdot \veps$, by \ref{9.2}(B).
Moreover, $\bkappa(T,3 \cdot \veps) \leq \bM(P) < k \cdot \veps$, according to \ref{8.0}(D) and \ref{9.2}(A).
\end{proof}

\begin{Corollary}
\label{9.6}
Assume that:
\begin{enumerate}
\item[(A)] $\Gamma > 0$ and $\kappa : \R^+_0 \to \R^+$ is non-increasing;
\item[(B)] $\calT \subset \calG(\bC^n)$ is $(\Gamma,\kappa)$-bounded;
\item[(C)] $0 < \veps \leq 1$;
\item[(D)] $k$ is a positive integer such that
\begin{equation*}
 \bc_{\ref{9.5}}(n) \cdot \left( \Gamma + \kappa ( \veps) \right) < k \cdot \veps.
\end{equation*}
\end{enumerate}
Then
\begin{equation*}
(\forall T \in \calT)(\exists P \in \calP^n_{k,\veps}) : \bG(T - P) < 3 \cdot \veps
\end{equation*}
\end{Corollary}

\section{Upper bound for the cardinality of $\calP^n_{k,\veps}$}

\begin{Empty}[Reduction to counting]
\label{10.1}
Recall that $\calP^n_{k,\veps}$ is defined in \ref{9.3}.
Note that $P \in \calP^n_{k,\veps}$ is entirely determined by the corresponding function $\theta : \Sigma^{n,0}_k \to Z(n,k,\veps)$.
This is equivalent to the knowledge of the function $f = \hat{\veps}(n,k,\veps)^{-1} \cdot \theta : \Sigma^{n,0}_k \to \Z$.
The latter satisfies two properties, namely
\begin{equation*}
\sum_{x \in \Sigma^{n,0}_k} f(x) = 0,
\end{equation*}
which corresponds to $\chi(P)=0$, and
\begin{equation*}
\sum_{x \in \Sigma^{n,0}_k} |f(x)| \leq k \cdot 2n \cdot (2k+1)^n,
\end{equation*}
which corresponds to $\bM(P) \leq k \cdot \veps$.
Therefore, $\calP^n_{k,\veps}$ is equipotent to the set
\begin{equation*}
\Z^{\Sigma^{n,0}_k} \cap \left\{ f : \sum f = 0 \text{ and } \sum |f| \leq k \cdot 2n \cdot (2k+1)^n \right\}.
\end{equation*}
In particular, $\rmcard \calP^n_{k,\veps}$ is independent of $\veps$.
\end{Empty}

\begin{Empty}[Counting]
\label{10.2}
Here, we are given two integers $1 \leq q \leq p$ and we consider the set
\begin{equation*}
E(p,q) = \Z^{\{1,\cdots,q\}} \cap \left\{ f : \sum_{j=1}^q f(j) = 0 \text{ and } \sum_{j=1}^q |f(j)| \leq p \right\}
\end{equation*}
as well as its cardinality
\begin{equation*}
\bc_{\theTheorem}(p,q) = \rmcard E(p,q).
\end{equation*}
Our goal is to find asymptotics for $\bc_{\theTheorem}(p,q)$ as $p,q \to \infty$ under the additional assumption that $q = o(p)$.
Indeed, we will apply this to $\bc_{\theTheorem}(p_k,q_k)$ as $k \to \infty$, where $q_k = (2k+1)^n$ and $p_k = k \cdot 2n \cdot (2k+1)^n$ as in \ref{10.1}.
\begin{enumerate}
\item[(A)] {\it The set of those $f : \{1,\cdots,q\} \to \Z$ such that $f(j) \geq 0$ for all $j$ and $\sum_{j=1}^q f(j) = p$ has cardinality $\bin{p+q-1}{q-1}$.}
\end{enumerate}
\par 
This is a standard exercise using the combinatorial ``balls and bars'' technique. 
The number of ways to put $p$ undistinguishable balls into $q$ drawers is $(p+q-1)!$ divided by $p!$ (the number of ways to enumerate the balls) and by $(q-1)!$ (the number of ways to enumerate the bars separating the drawers).\cqfd
\begin{enumerate}
\item[(B)] {\it The set of those $f : \{1,\cdots,q\} \to \Z$ such that $f(j) > 0$ for all $j$ and $\sum_{j=1}^q f(j) = p$ has cardinality $\bin{p-1}{q-1}$.}
\end{enumerate}
\par 
Recall that we assume $q \leq p$.
With each $f$ as in (B) we associate $g : \{1,\cdots,q\} \to \Z$ by means of the formula $g(j) = f(j) - 1$ for all $j$.
This correspondence establishes a bijection between these $f$ here and those $g$ in (A) with $p$ replaced by $p-q$.\cqfd 
\begin{enumerate}
\item[(C)] {\it $\displaystyle \bc_{\theTheorem}(p,q) = 1 + \sum_{r=1}^{q-1} \bin{q}{r} \sum_{\substack{s \geq r \\ 2s \leq p}}\bin{s-1}{r-1}\bin{q-r+s-1}{s}$.}
\end{enumerate}
\par 
The case when $f$ vanishes identically is counted in the first term.
Abbreviate $S = \{1,\cdots,q\}$.
From now on, we assume that $f$ is not identically zero so that $S^+ = \{f > 0\}$ is a non-empty proper subset of $S$ (recall that $\sum_S f = 0$).
Thus, upon letting $r = \rmcard S^+$, we have $1 \leq r \leq q-1$.
Given such $r$, there are $\bin{q}{r}$ choices for $S^+$.
Once the choice of $S^+$ is made, upon letting $s = \sum_{S^+} f$, we have $r \leq s$ and $2s \leq p$ (since $s = \sum_{S^+} f = - \sum_{S \setminus S^+} f$ and $p \geq \sum_S |f| = \sum_{S^+} f - \sum_{S \setminus S^+} f = s - (-s)$).
For $f|_{S^+}$, there are $\bin{s-1}{r-1}$ choices, according to (B).
It remains to choose $f|_{S \setminus S^+}$.
The number of possible such choices (for $-f$) is given by (A) applied with $s$ in place of $p$ and $q-r$ in place of $q$.\cqfd
\par 
We will not use (C) in the remaining part of this paper.
Instead we will obtain an upper bound for $\bc_{\theTheorem}(p,q)$.
\begin{enumerate}
\item[(D)] {\it The set of those $f : \{1,\cdots,q\} \to \Z$ such that $f(j) \geq 0$ for all $j$ and $\sum_{j=1}^q f(j) \leq p$ has cardinality $\bin{p+q}{q}$.}
\end{enumerate}
\par 
With such $f$ we associate $g : \{1,\cdots,q,q+1\}$ so that $g|_{\{1,\cdots,q\}}=f$ and $g(q+1) = p - \sum f$.
This establishes a bijection between between these $f$ here and those $g$ in (A) corresponding to the same $p$ and $q$ replaced by $q+1$.\cqfd
\begin{enumerate}
\item[(E)] {\it The set of those $f : \{1,\cdots,q\} \to \Z$ such that $\sum_{j=1}^q |f(j)| \leq p$ has cardinality bounded above by $2^q\bin{p+q}{q}$.}
\end{enumerate}
\par 
Let $\calE$ be the set of those $f$ considered here and $\calD$ the set of those $f$ considered in (D).
Then the map $\{-1,1\}^q \times \calD \to \calE$ that sends $(\la \veps_j \ra_{j=1}^q , f)$ to $j \mapsto \veps_j f(j)$ is surjective.\cqfd
\begin{enumerate}
\item[(F)] {\it $\displaystyle \bc_{\theTheorem}(p,q) \leq \frac{2^q(p+q)^{q}}{q!}$.}
\end{enumerate}
\par 
This follows the observation that $E(p,q)$ is contained in the set of those $f$ described in (E) and the trivial upper bound $\displaystyle \bin{n}{k} \leq \frac{n^k}{k!}$.\cqfd
\begin{enumerate}
\item[(G)] {\it $\displaystyle \ln \left( \bc_{\theTheorem}(p,q) \right) \leq q \cdot \ln \left( \frac{11 \cdot p}{q}\right)$.}
\end{enumerate}
\par 
Recalling that $q \leq p$ we infer from (F) that $\bc_{\theTheorem}(p,q) \leq \frac{(4p)^{q}}{q!}$.
Since $\ln(q!) \geq q \ln(q) - q +1$, we have
\begin{equation*}
\ln \left( \bc_{\theTheorem}(p,q) \right) \leq q \cdot  \ln(4p) - \ln(q!) \leq q \cdot [ \ln(4p) - \ln(q) +1]  = q \cdot  \ln \left( \frac{4 \cdot e \cdot p}{q} \right),
\end{equation*}
and we notice that $4\cdot e \leq 11$.\cqfd
\end{Empty}

\begin{Theorem}
\label{10.3}
For all positive integer $k$ and all $0 < \veps \leq 1$ we have
\begin{equation*}
\ln \left( \rmcard \calP^n_{k,\veps} \right) \leq (2k+1)^n \ln (22 \cdot n \cdot k).
\end{equation*}
\end{Theorem}

\begin{proof}
As explained before, we shall apply \ref{10.2}(G) with $q_k = (2k+1)^n$ and $p_k = k \cdot 2n \cdot q_k$.
Indeed, it follows from \ref{10.1} that $\rmcard \calP^n_{k,\veps} = \rmcard E(p_k,q_k)$ where the latter is defined in \ref{10.2}.
It then follows from \ref{10.2}(G) that 
\begin{multline*}
\ln \left( \rmcard \calP^n_{k,\veps} \right) = \ln \bc_{\ref{10.2}}(p_k,q_k) \leq q_k \cdot \ln \left( \frac{11 \cdot p_k}{q_k }\right) \leq q_k \cdot \ln \left( \frac{11 \cdot p_k}{q_k}\right).
\end{multline*}
\end{proof}

\begin{Remark}
The inequality $\rmcard \calP^n_{k,\veps} \leq (22 \cdot n \cdot k)^{(2k+1)^n}$ is a slight improvement on the brute force estimate $\rmcard \calP^n_{k,\veps} \leq \left[ 1 + 4n \cdot k \cdot (2k+1)^n \right]^{(2k+1)^n}$ that we derive from the fact that (in the notation of \ref{10.2}) $E(p,q) \subset \{-p,\cdots, 0, \cdots p\}^{\{1,\cdots,q\}}$.
The improvement is due to the condition $\bM(P) \leq k \cdot \veps$ and owes nothing to the condition $\chi(P)=0$, both from the definition of $\calP^n_{k,\veps}$.
\end{Remark}

\begin{Corollary}
\label{10.4}
Assume that:
\begin{enumerate}
\item[(A)] $\Gamma > 0$ and $\kappa : \R^+_0 \to \R^+$ is non-increasing;
\item[(B)] $\calT \subset \calG(\bC^n)$ is $(\Gamma,\kappa)$-bounded;
\item[(C)] $0 < \veps \leq 1$;
\item[(D)] $k$ is a positive integer such that
\begin{equation*}
 \bc_{\ref{9.5}}(n) \cdot \left( \Gamma + \kappa \left( \frac{\veps}{6} \right) \right) < k \cdot \left( \frac{\veps}{6} \right);
\end{equation*}
\item[(E)] $\calQ \subset \calT$ is an $\veps$-net, \ie
\begin{equation*}
(\forall Q_1 \in \calQ)(\forall Q_2 \in \calQ): Q_1 \neq Q_2 \Rightarrow \bG(Q_1-Q_2) \geq \veps.
\end{equation*}
\end{enumerate}
Then
\begin{equation*}
\rmcard \calQ \leq [22 \cdot n \cdot k ]^{(2k+1)^n}.
\end{equation*}
\end{Corollary}

\begin{proof}
It follows from \ref{9.6} that each $T \in \calT$ can be associated with some $P \in \calP^n_{k,\frac{\veps}{6}}$ such that $\bG(T-P) < \frac{\veps}{2}$.
As $\calQ \subset \calT$, this defines a map $\calQ \to \calP^n_{k,\frac{\veps}{6}}$.
This map is injective: If $Q_1$ and $Q_2$ are associated with the same $P$ then $\bG(Q_1 - Q_2) \leq \bG(Q_1-P) + \bG(Q_2-P) < \veps$, therefore, $Q_1=Q_2$.
The conclusion is now a consequence of \ref{10.3}.
\end{proof}


\bibliographystyle{amsplain}
\bibliography{/Users/thierry/Documents/LaTeX/Bibliography/thdp}




\end{document}